\documentclass[reqno]{amsart}
 \usepackage{amsbsy,amssymb,amscd,amsfonts,latexsym,amstext,delarray,
 amsmath,graphicx, color}
 \usepackage[dvips]{epsfig}
\usepackage{hyperref}
\input xypic

\definecolor{Green}{rgb}{0,1,0}
\definecolor{Blue}{RGB}{0,0,191}
\definecolor{mathmodecolor}{RGB}{0,102,0}
\definecolor{keywordcolor}{RGB}{0,51,151}
\definecolor{sourcebackgroundcolor}{RGB}{255,247,223}
\definecolor{unixagred}{RGB}{255,0,0}
\definecolor{lightgray}{RGB}{191,191,191}
\definecolor{green}{RGB}{1,191,191}

\newtheorem{thm}{Theorem}[section]
\newtheorem{prop}[thm]{Proposition}
\newtheorem{cor}[thm]{Corollary}
\newtheorem{lem}[thm]{Lemma}

\newtheorem{defn}[thm]{Definition}
\newtheorem{rem}[thm]{Remark}

\def\qqq{\,,\quad~\forall}

\def\Aut{{\rm Aut}}

\def\GL{{\rm GL}}

\def\Hom{{\rm Hom}}

\def\\adsod{{\rm \adsod}}

\def\Res{{\rm Res}}
\def\sign{{\rm sign}}

\def\Spec{{\rm Spec\,}}
\def\Sp{{\rm Spec}\,}

\def\Tr{{\rm Tr}}

\def\A{{\mathbb A}}

\def\C{{\mathbb C}}
\def\F{{\mathbb F}}
\def\K{{\mathbb K}}

\def\N{{\mathbb N}}
\def\P{{\mathbb P}}
\def\Q{{\mathbb Q}}
\def\R{{\mathbb R}}
\def\Z{{\mathbb Z}}
\def\H{{\mathbb H}}

\def\Tr{{\rm Tr}}

\def\cD{{\mathcal D}}

\def\cH{{\mathcal H}}

\def\cL{{\mathcal L}}
\def\c\ads{{\mathcal \ads}}

\def\cO{{\mathcal O}}
\def\cP{{\mathcal P}}

\def\cS{{\mathcal S}}

\def\cU{{\mathcal U}}

\def\cW{{\mathcal W}}

\newcommand{\ie}{{\it i.e.\/}\ }

\newcommand{\cf}{{\it cf.\/}\ }

\newcommand{\resp}{{\it resp.\/}\ }

\def\dim{{\mbox{dim}}}

\def\Hom {{\mbox{Hom}}}

\def\ffp{\mathfrak{p}}

\def\strong{{\cS_\infty}}

\def\\adso{\mathfrak{\adso}}
\def\An{\mathfrak{ Ring}}
\def\han{\mathfrak{ Hring}}

\def\\adsr{\mathfrak{ \adsR}}

\def\Se{\mathfrak{ Sets}}

\def\ss{{\rm sign}\,}

\def\urep{\vartheta}
\def\rep{\vartheta}
\def\spad{{\P^1_{\F_1}}}
\def\spadu{{\P^1_{\F_1}}}
\def\Omm{\Omega}

\def\rmax{\R_+^{\rm max}}
\def\ads{\H_\K}
\def\kras{\mathbf{K}}
\def\sign{\mathbf{S}}
\def\Mo{\mathfrak{Mo}}
\def\rco{\R^{\rm{convex}}}
\def\ad{ad\`ele }

\newcommand{\nil}[1]{}

\title
{From monoids to hyperstructures: in search
of an absolute arithmetic.}

\author[Connes]{Alain Connes}
\author[Consani]{Caterina Consani}
\address{A.~Connes: Coll\`ege de France \\
3, rue d'Ulm \\ Paris, F-75005 France
\\ I.H.E.S. and Vanderbilt
University} \email{alain\@@connes.org}
\address{C.~Consani: Mathematics Department \\ Johns Hopkins
University \\ Baltimore, MD 21218 USA} \email{kc\@@math.jhu.edu}
\thanks{The authors are partially supported by the NSF grant
DMS-FRG-0652164.}


\begin{document}

\maketitle

\begin{abstract} We show that the trace formula interpretation of the explicit formulas
 expresses the counting function $N(q)$ of the hypothetical curve $C$ associated to the Riemann zeta function, as an intersection number involving
the scaling action on the \ad class space. Then, we discuss the algebraic structure of the \ad class space both as a monoid and as a hyperring. We construct an extension $\rco$ of the hyperfield $\sign$ of signs, which is the hyperfield analogue of the semifield  $\R_+^{\rm max}$ of tropical geometry, admitting a one parameter group of automorphisms fixing $\sign$. Finally, we develop  function theory over $\Spec\kras$ and we show how to recover the field of real numbers from a purely algebraic construction, as the function theory over $\Spec\sign$.
\end{abstract}

\tableofcontents

\section{Introduction}

The geometry of  algebraic curves underlying the structure of global fields of positive characteristic has shown its crucial role in the process of solving several fundamental questions in number theory. Furthermore, some combinatorial formulas, such as the equation supplying the cardinality of the set of rational points of a Grassmannian  over a finite field $\F_q$, are known to be rational expressions keeping a meaningful and interesting value also when $q=1$.
These results motivate the search for a mathematical object which is expected to be a non-trivial limit of  Galois fields  $\F_q$, for $q=1$. The goal is to define an analogue, for number fields, of the geometry underlying the arithmetic theory of the function fields.
 Inspired by this problem and by the pioneering work of  R. Steinberg and J. Tits, C. Soul\'e  has associated a zeta function to any sufficiently regular counting-type function $N(q)$, by considering the following limit
\begin{equation}\label{zetadefn}
\zeta_N(s):=\lim_{q\to 1}Z(q,q^{-s}) (q-1)^{N(1)}\qquad s\in\R.
\end{equation}
Here,  $Z(q,q^{-s})$ denotes the evaluation, at $T=q^{-s}$,  of the  Hasse-Weil zeta function
\begin{equation}\label{zetadefn1}
Z(q,T) = \exp\left(\sum_{r\ge 1}N(q^r)\frac{T^r}{r}\right).
\end{equation}
For the consistency of the formula \eqref{zetadefn}, one requires that the counting function $N(q)$ is defined for all real numbers
 $q \geq 1$ and not only for prime integers powers as for the counting function in \eqref{zetadefn1}.  For many simple examples of algebraic varieties, like the projective spaces, the function  $N(q)$ is known to extend unambiguously  to all real positive numbers. The associated zeta function $\zeta_N(s)$ is easy to compute and it produces the expected outcome. For a projective line, for example, one finds  $\zeta_N(s) = \frac{1}{s(s-1)}$. Another case which is easy to handle and also carries a richer geometric structure is provided by a Chevalley group scheme. The study of these varieties has shown that these schemes are endowed with a basic (combinatorial) skeleton structure and, in   \cite{ak}, we proved that they are varieties over $\F_1$, as defined by Soul\'e  in \cite{Soule}.
 
 To proceed further, it is natural to wonder on the existence of a  suitably defined curve $C=\overline {\Sp\,\Z}$  over $\F_1$, whose zeta function $\zeta_C(s)$ agrees with the {\em complete} Riemann zeta function  $\zeta_\Q(s)=\pi^{-s/2}\Gamma(s/2)\zeta(s)$ (\cf also \cite{Ku} and \cite{Manin}). Following the interpretation of $N(1)$ given in \cite{Soule}, this value should coincide with the Euler characteristic
of the curve $C$. However, since  $C$ should have infinite genus, one deduces that $N(1)=-\infty$, in apparent contradiction with the expected positivity of $N(q)$, for $q>1$. This result also prevents one to use the limit definition \eqref{zetadefn}. In \cite{announc3}, we have shown how to solve these difficulties by replacing the limit definition \eqref{zetadefn} with an integral formula and by computing explicitly the {\em distribution} $N(q)$ which fulfills the expected positivity, for $q>1$, and the divergence at $q=1$.  

In section \ref{countoad} of this paper,  we show
 how to implement the trace formula understanding of the explicit formulas in number-theory,
to express the counting function $N(q)$ as an {\em intersection number} involving
the scaling action of the id\`ele class group on the \ad class space.  This description is the natural corollary of  an earlier result proved by the first author of this paper  in \cite{Co-zeta} (and motivated  by \cite{gui}), showing that the Riemann-Weil explicit formulas inherit a trace formula meaning when one  develops analysis on the noncommutative space of the \ad classes (we also refer to \cite{CCM}, \cite{CMbook} and \cite{Meyer} for further work on the same topic).\vspace{.02in}

In \cite{Co-zeta}, as well as in the above papers, the \ad class space $\A_\K/\K^\times$ of a global field $\K$ has been studied as a noncommutative space. Only very recently (\cf\cite{wagner}),  we have been successful to associate an algebraic structure to $\A_\K/\K^\times$ using which this space  finally reveals its deeper nature of a {\em hyperring of functions}, likewise the space-time geometry  in quantum field theory which manifests itself through  functional analysis.
The hyperring structure of $\A_\K/\K^\times$ has fully emerged by combining the following two properties:\vspace{.05in}

$-$~$\A_\K/\K^\times$ is a commutative monoid, so that one may apply to this space the geometry of monoids of K.~ Kato and A.~Deitmar.\vspace{.05in}

$-$~$\A_\K/\K^\times$ is a hyperring over the Krasner hyperfield $\kras=\{0,1\}$.\vspace{.05in}

 In section \ref{sectmonoids}, we describe the first structure associated to this space.
In particular, we overview  several of our recent results (\cf \cite{announc3}) showing that the natural monoidal structure on the ad\`ele class space, when combined with
one of the simplest examples of  monoidal schemes \ie the projective line $\spad$,
 provides a  geometric framework to understand conceptually the spectral realization of the zeros of $L$-functions, the functional equation and the explicit formulas. All these highly non-trivial statements appear naturally by simply computing the cohomology of a natural sheaf $\Omm$ of functions on the set   of rational points of the monoidal scheme $\spad$ on  the monoid $M=\A_\K/\K^\times$  of \ad classes.  The  sheaf  $\Omm$ is a sheaf of complex vector spaces over the geometric realization of the functor associated to the projective line.
It is a striking fact  that despite the apparent simplicity of the construction of $\spad$ the computation of $H^0(\spad,\Omm)$, already yields the graph of the Fourier transform. The cohomology
    $H^0(\spad,\Omm)$ is given, up to a finite dimensional space, by the graph  of the Fourier transform
    acting on the co-invariants for the action of $\K^\times$ on the  Bruhat-Schwartz space $\cS(\A_\K)$. Moreover, the spectrum of the natural action   of the id\`ele class group $C_\K$ on the cohomology $\displaystyle{H^1(\spad,\Omm)}$ provides the spectral realization of the zeros of Hecke  $L$-functions.

In section \ref{hyper}, we review the hyperring structure associated to the ad\`ele class space.  In \cite{wagner}, we proved that the \ad class space possesses  a rich additive structure which plays an essential role to provide the correct arithmetic setup on this space and to obtain, in positive characteristic, a canonical identification of the groupoid of prime elements of the hyperring $\A_\K/\K^\times$ with the loop groupoid of the maximal abelian cover of the algebraic curve underlying the arithmetic of the function field $\K$. It is an interesting coincidence that the first summary on hyperring theory, due to M. Krasner (\cf\cite{Krasner}), appeared in the same proceeding volume together with the seminal paper of J. Tits \cite{Tits} where he introduced ``le corps de caract\'eristique un". The distinction between the algebraic structure that Tits proposed as the degenerate case of $\F_q$ for $q=1$, \ie  ``le corps form\'e du seul \'el\'ement   $1=0$",  and the Krasner hyperfield $\kras=\{0,1\}$ is simply that in $\kras$ one keeps the distinction $1\neq 0$, while recording the degeneracy  by allowing the sum $1+1$ to be maximally ambiguous. 

In section \ref{hypers}, we recall some of our results contained in \cite{wagner} showing, for instance, that in spite of the apparent naivety  of the definition of the Krasner hyperfield $\kras$, the classification of its finite extensions is equivalent to a deep open question in classical projective geometry.

 When $\K$ is a number field, there is a basic main question that needs to be addressed, namely the search for a  substitute, in characteristic zero, of the algebraic closure $\bar\F_q$ of the field of constants  of a function field and the understanding of its role in the geometric construction of the curve associated to that field.
 In section \ref{rconvex}, we show that the hyperfield $\sign=\{-1,0,1\}$ of signs admits an infinite hyperfield extension $\rco$ which is obtained  from $\R$ by suitably altering the rule of addition. This extension has characteristic one (\ie $x+x=x$ for all $x\in \rco$) and it contains $\sign$ as the smallest sub-hyperfield. The group of automorphisms $\Aut(\rco)$ is the multiplicative group of non-zero  real numbers acting on the hyperfield $\rco$ by exponentiation. This group plays the role of the cyclic group generated by the Frobenius automorphism acting on $\bar\F_p$, in characteristic $p>1$.\vspace{.02in}

 In section \ref{fnthhr}, we develop the generalities of the study of the function theory on the affine spectrum $\Spec(R)$ of a hyperring $R$. The novelty (and subtlety) with respect to the classical case of (commutative) rings arises from the failure in this framework, of the isomorphism (holding for rings)
 $$
 \Hom(\Z[T],R)\simeq R.
 $$
In the classical case, the above isomorphism provides one with an identification of the elements of $R$ with functions on $\Sp(R)$, understood as  morphisms to the affine line $\cD=\Sp(\Z[T])$. In the context of hyperrings, we define the functions on $\Sp(R)$ as the elements of the set $\cD(R)=\Hom(\Z[T],R)$. We  implement the natural coproduct structures on the Hopf algebra $\cH=\Z[T]$ (corresponding to the addition: $\Delta^+(T)=T\otimes 1+1\otimes T$ and the multiplication $\Delta^\times(T)=T\otimes T$) to define operations on functions on $\Sp(R)$.
In sections \ref{sectspeckras} and \ref{sectspecsign}, we investigate the outcoming hyperstructures on  $\cD(\kras)$ and $\cD(\sign)$. The natural morphism from $\cD$ to $\Spec\,\Z$ determines a restriction map $\pi:\cD(\kras)=\Spec(\Z[T])\to \Hom(\Z,\kras)=\Spec(\Z)$. Then, one sees that the two hyperoperations of addition and multiplication on $\cD(\kras)$ take place fiber-wise \ie within each fiber $\pi^{-1}(p)$, $p\in \Spec(\Z)$. Theorem \ref{thmbijp} asserts that, for a {\em finite and proper} prime  $p$,
 the hyperstructure on the fiber $\pi^{-1}(p)$  coincides with the quotient hyperring structure on the orbit-set $\Omega/\Aut(\Omega)$, where $\Omega$ is an algebraic closure of the field of fractions $\F_p(T)$.
 Theorem \ref{thmalg} states that the fiber over the {\em generic} point of $\Spec(\Z)$ contains the quotient hyperstructure $\bar\Q/\Aut(\bar\Q)$ where $\bar \Q$ is the algebraic closure of $\Q$,  although the operations involving the generic point are more subtle: \cf Theorems \ref{genericplus} and \ref{genericprod}.
 
 We expect that a similar development will hold when the Hopf algebra $\Z[T]$  is replaced by the Hopf algebra of a Chevalley group.\vspace{.02in}

In  section \ref{sectspecsign}, we show that the hyperstructure $\cD(\sign)$ defines a slight refinement of the field $\R$ of real numbers. To each element $\varphi\in \cD(\sign)$ corresponds a real number $Re(\varphi)\in[-\infty,\infty]$ given as a Dedekind cut. The map $Re:  \cD_{\rm finite}(\sign)\to \R$ is a surjective homomorphism
whose kernel is an ideal isomorphic to $\sign$. A slight difference between the real numbers and  the subset  $\cD_{\rm finite}(\sign)$ of the finite elements of $\cD(\sign)$  only occurs over real algebraic numbers which have three representatives given by the three homomorphisms:
$$
\Z[T]\to\sign,\qquad P(T)\mapsto \lim_{\epsilon\to 0+}\ss P(\alpha+t\epsilon)\,, \ t\in\{-1,0,1\}.
$$
The richness and complexity  of the structure of functions on hyperring spectra,
even in the simplest examples of $\Sp(\kras)$ and $\Sp(\sign)$, together with the
construction of the hyperfield $\rco$, are both clear indications of the existence of an interesting and yet unexplored
 arithmetic theory which has just started to emerge from the basic principles outlined in this work.

\section{From the counting function to the ad\`ele class space}\label{countoad}

In this section we show that the trace formula interpretation of the explicit formulas
 expresses the counting function $N(q)$, of the hypothetical curve $C$ associated to the complete Riemann zeta function, as an intersection number involving
the scaling action of the id\`ele class group on the \ad class space.

 \subsection{The counting function of $C=\overline {\Sp\,\Z}$} \label{distrsubsect}

 As explained in the introduction, it is natural to wonder on the existence of a  suitably defined ``curve'' $C=\overline {\Sp\,\Z}$  over $\F_1$, whose zeta function $\zeta_C(s)$ is the {\em complete} Riemann zeta function  $\zeta_\Q(s)=\pi^{-s/2}\Gamma(s/2)\zeta(s)$ (\cf also \cite{Manin}). To by-pass the difficulty inherent to the definition \eqref{zetadefn}, when $N(1)=-\infty$,  one works  with the logarithmic derivative
\begin{equation}\label{normael}
    \frac{\partial_s\zeta_N(s)}{\zeta_N(s)}=-\lim_{q\to 1} F(q,s)
\end{equation}
where
\begin{equation}\label{fqsdefn}
F(q,s)=-\partial_s \sum_{r\ge 1}N(q^r)\frac{q^{-rs}}{r}\,.
\end{equation}
Then one finds (\cf\cite{announc3} Lemma 2.1), under suitable regularity conditions on $N(u)$, that

\begin{lem} \label{compute1}With the above notations and for $\Re\mathfrak e(s)$ large enough, one has
\begin{equation}\label{lim}
    \lim_{q\to 1} F(q,s) = \int_1^\infty N(u)u^{-s}d^*u\,,\ \ d^*u=du/u
\end{equation}
and
\begin{equation}\label{logzetabis}
    \frac{\partial_s\zeta_N(s)}{\zeta_N(s)}=-\int_1^\infty  N(u)\, u^{-s}d^*u\,.
\end{equation}
\end{lem}

 The integral equation \eqref{logzetabis} is more manageable and general  than the limit formula \eqref{zetadefn}. It produces a precise description of the counting function $N_C(q)=N(q)$ associated to $C$. In fact, \eqref{logzetabis} shows in this case that
\begin{equation}\label{special}
   \frac{\partial_s\zeta_\Q(s)}{\zeta_\Q(s)}=-\int_1^\infty  N(u)\, u^{-s}d^*u\,.
\end{equation}
To determine explicitly $N(u)$, one uses  the Euler product for $\zeta_\Q(s)$ and when $\Re\mathfrak e(s)>1$, one derives
\begin{equation}\label{special1}
   -\frac{\partial_s\zeta_\Q(s)}{\zeta_\Q(s)}=\sum_{n=1}^\infty\Lambda(n)n^{-s}+\int_1^\infty  \kappa(u)\, u^{-s}d^*u\,.
\end{equation}
 Here, $\Lambda(n)$ is the von-Mangoldt function taking the value $\log p$ at prime powers $p^\ell$  and zero otherwise.  $\kappa(u)$ is the distribution
defined, for any test function $f$, as
\begin{equation}\label{kappadu}
\int_1^\infty\kappa(u)f(u)d^*u=\int_1^\infty\frac{u^2f(u)-f(1)}{u^2-1}d^*u+cf(1)\,, \qquad c=\frac12(\log\pi+\gamma)
\end{equation}
where $\gamma=-\Gamma'(1)$ is the Euler constant. The distribution $\kappa(u)$ is positive on $(1,\infty)$ where, by construction, it is given by $\kappa(u)=\frac{u^2}{u^2-1}$.
Hence, we deduce that the counting function $N(q)$ of the hypothetical curve $C$ over $\F_1$, is the {\em distribution} defined by the sum of $\kappa(q)$ and a discrete term given by the derivative taken in the sense of distributions, of the function\footnote{the value at the points of discontinuity does not affect the distribution}
\begin{equation}\label{varphifunc}
    \varphi(u)=\sum_{n<u}n\,\Lambda(n).
\end{equation}
Indeed, since  $d^*u=\frac{du}{u}$, one has for any test function $f$,
$$
\int_1^\infty f(u)\left(\frac{d}{du}\varphi(u)\right)d^*u=\int_1^\infty \frac{f(u)}{u}d\varphi(u)=\sum \Lambda(n)f(n).
$$
Thus one can write \eqref{special1} as
\begin{equation}\label{special2}
   -\frac{\partial_s\zeta_\Q(s)}{\zeta_\Q(s)}=\int_1^\infty \left(\frac{d}{du}\varphi(u)+ \kappa(u)\right)\, u^{-s}d^*u.
\end{equation}
If one compares the equations \eqref{special2} and \eqref{special}, one derives the following formula for $N(u)$
\begin{equation}\label{Nu}
    N(u)=\frac{d}{du}\varphi(u)+ \kappa(u).
\end{equation}
One can then use the explicit formulas to express $\varphi(u)$ in terms of  the set $Z$ of non-trivial zeros of the Riemann zeta function. One has the formula (\cf~\cite{Ingham}, Chapter IV, Theorems 28 and 29) valid for $u>1$ (and not a prime power)
\begin{equation}\label{ingham1}
\varphi(u)=\frac{u^2}{2}-\sum_{\rho\in Z}{\rm order}(\rho)\frac{u^{\rho+1}}{\rho+1}+a(u).
\end{equation}
Here, one sets
\begin{equation}\label{inghamb}
a(u)={\rm ArcTanh}(\frac 1u)- \frac{\zeta'(-1)}{\zeta(-1)}.
\end{equation}
 Notice that  the sum over $Z$ in \eqref{ingham1} has to be taken in a symmetric manner to ensure the convergence, \ie as a limit of the partial sums over the  symmetric set $Z_m$ of first $2m$ zeros. When one differentiates \eqref{ingham1} in a formal way, the term in $a(u)$ gives
$$
\frac{d}{du}a(u)=\frac{1}{1-u^2}.
$$
Hence, at the formal level \ie by disregarding the principal value, one obtains
$$
\frac{d}{du}a(u)+\kappa(u)=1.
$$
Thus, after a formal differentiation of \eqref{ingham1}, one deduces
\begin{equation}\label{formallevel}
    N(u)=\frac{d}{du}\varphi(u)+ \kappa(u)\sim u-\sum_{\rho\in Z}{\rm order}(\rho)\,u^{\rho}+1.
\end{equation}

Notice that in the above formal computations we have neglected to consider the principal value for the distribution $\kappa(u)$. By taking this into account, we obtain the following more precise result (for the proof we refer to \cite{announc3}, Theorem 2.2)
\begin{thm}\label{fine1}
The tempered distribution $N(u)$ satisfying the equation
$$
-\frac{\partial_s\zeta_\Q(s)}{\zeta_\Q(s)}=\int_1^\infty  N(u)\, u^{-s}d^*u\,
$$
is positive on $(1,\infty)$ and on $[1,\infty)$ is given  by
\begin{equation}\label{fin2}
    N(u)=u-\frac{d}{du}\left(\sum_{\rho\in Z}{\rm order}(\rho)\frac{u^{\rho+1}}{\rho+1}\right)+1
\end{equation}
where the derivative is taken in the sense of distributions, and the value at $u=1$ of the  term
 $\displaystyle{\omega(u)=\sum_{\rho\in Z}{\rm order}(\rho)\frac{u^{\rho+1}}{\rho+1}}$ is given  by
\begin{equation}\label{definethet}
\omega(1)=\sum_{\rho\in Z}{\rm order}(\rho)\frac{1}{\rho+1}=\frac 12+ \frac \gamma 2+\frac{\log4\pi}{2}-\frac{\zeta'(-1)}{\zeta(-1)}.
\end{equation}
\end{thm}

This result supplies a strong indication on the coherence of the quest for an arithmetic theory over $\F_1$. For an irreducible, smooth and projective algebraic curve $X$ over a prime field $\F_p$, the counting function is of the form
 $$
\#X(\F_q)=N(q)=q-\sum_{\alpha} \alpha^r+1,\qquad   q=p^r
$$
where the numbers $\alpha$'s are the complex roots  of the characteristic polynomial of the Frobenius endomorphism acting on the \'etale cohomology $H^1(X\otimes\bar\F_p,\Q_\ell)$ of the curve ($\ell\neq p$).
By writing these roots in the form $\alpha=p^\rho$, for $\rho$ a zero of the Hasse-Weil zeta function of $X$, the above equality reads as
\begin{equation}\label{inghamnsimilar}
\#X(\F_q)=N(q)=q- \sum_{\rho}{\rm order}(\rho)\, q^{\rho} +1.
\end{equation}
 The equations \eqref{fin2} and \eqref{inghamnsimilar} are now completely identical, except for the fact that in \eqref{inghamnsimilar} the values of $q$ are restricted to the discrete set of powers of $p$ and that \eqref{inghamnsimilar} involves only a finite sum, which allows one to differentiate term by term.

\subsection{Explicit formulas} \label{explicit}

Equation \eqref{ingham1} is a typical application of the Riemann-Weil explicit formulas. These formulas become natural when lifted to the id\`ele class group. In this section we show that, even if a definition of the hypothetical curve $C$ is at this time still out of reach, its counterpart, through the application of the class-field theory isomorphism, can be realized by a space of adelic nature and  in agreement with  some earlier constructions of Connes, Consani and Marcolli: \cf~\cite{CCM}, \cite{CCM2} and \cite{ccm}.\vspace{.05in}

We start by considering the explicit formulas in the following concrete form (\cf~\cite{Cartier}).
Let $F(u)$ be a function defined on $[1,\infty)$ and such that $F(u)=O(u^{-1/2-\epsilon})$. Then one sets
\begin{equation}\label{expl}
\Phi(s)=\int_1^\infty F(u)\,u^{s-1}du.
\end{equation}
The explicit formula  takes the following form
\begin{equation}\label{explfor}
\int_1^\infty (u^{-1/2}+u^{-3/2})F(u)du-\sum_{\rho\in Z}{\rm order}(\rho)\Phi(\rho-\frac 12)=
\end{equation}
$$
=\sum_p\sum_{m=1}^\infty \log p \,\,p^{-m/2}F(p^m)+(\frac \gamma 2+\frac{\log\pi}{2})F(1)
+\int_1^\infty\frac{t^{3/2}F(t)-F(1)}{t(t^2-1)}dt.
$$
We apply this formula with the function $F_x$ determined by the conditions
\begin{equation}\label{explfunc}
F_x(u)=u^{\frac 32} \ \  \text{for} \ u\in [1,x]\,, \ \ F_x(u)=0 \  \text{for} \ u>x.
\end{equation}
Then, we obtain
\begin{equation}\label{explphi}
\Phi_x(s)=\int_1^\infty F_x(u)\,u^{s-1}du=\int_1^x u^{\frac 32}\,u^{s-1}du=\frac{x^{3/2+s}}{3/2+s}-\frac{1}{3/2+s}.
\end{equation}
Thus, it follows that
\begin{equation}\label{explphi1}
\Phi_x(\frac 12)=\frac{x^2}{ 2}-\frac 12\,, \ \
\Phi_x(-\frac 12)=x-1\,, \ \Phi_x(\rho-\frac 12)=
\frac{x^{1+\rho}}{1+\rho}-\frac{1}{1+\rho}.
\end{equation}
The left-hand side of the explicit formula \eqref{explfor} gives, up to a constant
\begin{equation}\label{expl2}
J(x)=\frac{x^2}{ 2}+x -\sum_{\rho\in Z}{\rm order}(\rho)\frac{x^{1+\rho}}{1+\rho}.
\end{equation}
The first term on the right-hand side of \eqref{explfor} gives
\begin{equation}\label{expl3}
  \varphi(x)=\sum_{n<x}n\,\Lambda(n)
  \end{equation}
  while the integral on the right-hand side of \eqref{explfor} gives
  \begin{equation}\label{expl4}
 \int_1^\infty\frac{t^{3/2}F(t)-F(1)}{t(t^2-1)}dt=x- {\rm ArcTanh}(\frac 1x)+ {\rm constant}.
\end{equation}
Thus the explicit formula \eqref{explfor} is transformed into the equality
\begin{equation}\label{expl5}
\sum_{n<x}n\,\Lambda(n)=\frac{x^2}{ 2} -\sum_{\rho\in Z}{\rm order}(\rho)\frac{x^{1+\rho}}{1+\rho}
+ {\rm ArcTanh}(\frac 1x)+ {\rm constant}.
\end{equation}
This formula is the same as \eqref{ingham1}. We refer to \cite{Ingham} for a precise justification of the analytic steps.
It follows that the left-hand side \eqref{expl2} of the explicit formula gives a natural primitive $J(x)$
of the counting function $N(x)$. It is thus natural to differentiate formally the family of functions $F_x$ with respect to $x$ and see what the right-hand side of the explicit formula is transformed into. By construction, one has, for $u\geq 1$
$$
F_x(u)=u^{\frac 32}Y(u-x)
$$
where $Y$ is the characteristic function of the interval $(-\infty, 0]$. The derivative of $Y(s)$ is $-\delta(s)$. Thus, at the formal level, one derives
$$
\partial_x F_x=u^{\frac 32}\delta(u-x).
$$
In fact, it is more convenient to rewrite the explicit formula \eqref{explfor} in terms of the function $g(u)=u^{-\frac 12}F(u)$.
One then lets
\begin{equation}\label{explg}
\hat g(s)=\int_1^\infty g(u)\,u^{s}d^*u.
\end{equation}
The explicit formula then takes the form
\begin{equation}\label{explforg}
\hat g(0)+\hat g(1)-\sum_{\rho\in Z}{\rm order}(\rho)\hat g(\rho)
\end{equation}
$$
=\sum_p\sum_{m=1}^\infty \log p \,\,g(p^m)+(\frac \gamma 2+\frac{\log\pi}{2})g(1)
+\int_1^\infty\frac{t^{2}g(t)-g(1)}{t^2-1}d^*t.
$$
The function $g_x(u)$ corresponding to $\partial_x F_x$ is just $g_x(u)=u\delta_x(u)$ and
it is characterized, as a distribution, by its evaluation on test functions $b(x)$. This gives
\begin{equation}\label{deltafun}
    \int b(u)g_x(u)d^*u=b(x).
\end{equation}
Next, we show how to implement the trace formula interpretation of the explicit formulas to describe the counting function $N(u)$ as an intersection number. First we notice that the above explicit formula is a special case of the Weil explicit formulas. One
lets  $\K$ be a global field, $ \alpha$ a nontrivial character of
 $\A_\K/\K$ and $ \alpha = \prod \,  \alpha_v$ its local factors. Let $h \in \cS (C_\K)$ have compact support. Then
\begin{equation}\label{weil4}
\hat h (0) + \hat h (1)  - \sum_{\chi\in \widehat{C_{\K,1}}}\,
\sum_{Z_{\tilde\chi}} \hat h (\tilde\chi , \rho) = \sum_v
\int'_{\K_v^*} \frac{h(u^{-1}) }{   \vert 1-u   \vert} \, d^* u
\end{equation}
where  $\int'$ is normalized by $ \alpha_v$ and
$
\hat h (\chi , z) = \int h(u) \,
\chi (u) \,   \vert u   \vert^z \, d^* u.
$
These formulas become a trace formula whose geometric side  is of the form
\begin{equation}\label{geomside}
\Tr_{\rm distr}\left(\int h(u)\urep(u)d^*u\right )=\sum_v\int_{\K^*_v}\,\frac{h(u^{-1})}{|1-u|}\,d^*u.
\end{equation}
Here $
\urep(u)\xi(x)=\xi(u^{-1}x)
$ is the scaling action of the id\`ele class group $C_\K$  on the \ad class space $M=\A_\K/\K^*$. The subgroups $\K^*_v\subset C_\K$ appear as isotropy groups. One can understand why the terms $\displaystyle \frac{h(u^{-1})}{|1-u|}$ occur in the trace formula by computing formally as follows the trace  of the scaling operator $T=\theta(u^{-1})$
$$T\xi(x)=\xi(u x)=\int k(x,y)\xi(y)dy\,,
 $$
 given by the distribution kernel $k(x,y)=\delta(u x-y)$,
 $$
\Tr_{\rm distr}(T)=\int k(x,x)\,dx=\int \delta(u x-x)\,dx=\frac{1}{|u-1|}\int \delta(z)\,dz=\frac{1}{|u-1|}\,.
$$
We refer to \cite{Co-zeta}, \cite{Meyer} and \cite{CMbook} for the detailed treatment.
 We apply \eqref{geomside}  by taking $\K=\Q$
  and the function $h$ of the form $h(u)=g(|u|)$ where the support of the function $g$ is contained in
$(1,\infty)$. On the left hand side of \eqref{geomside} one first performs the integration in the kernel $C_{\Q,1}$ of the module $C_{\Q}\to \R_+^*$. At the geometric level, this corresponds to taking the quotient of $M$ by the action of $C_{\Q,1}$. We denote by $\urep_u$ the scaling action on this quotient. By construction this action only depends upon $|u|\in \R_+^*$.
The equality \eqref{deltafun} means that when we consider the distributional trace of an expression
of the form
$
\int  g_x(u) \urep_u d^*u
$
we are in fact just taking the distributional trace of $\urep_x$ since
\begin{equation}\label{base}
    \int  g_x(u) \urep_u d^*u=\urep_x
\end{equation}
 thus we are simply considering
an intersection number.
We now look at the right hand side of \eqref{geomside}, \ie at the terms
\begin{equation}\label{weil2}
\int'_{\K_v^*} \frac{h(u^{-1}) }{   \vert 1-u   \vert} \, d^* u.
\end{equation}
Since $h(u)=g(|u|)$ and  the support of the function $g$ is contained in
$(1,\infty)$, one sees that the integral \eqref{weil2} can be restricted  in all cases to the unit
ball $\{u\,; \,|u|<1\}$ of the local field $\K_v$. In particular, for the finite places one has
 $\vert 1-u   \vert=1$, thus for each finite prime $p\in\Z$ one has
\begin{equation}\label{weil3}
\int'_{\Q_p^*} \frac{h(u^{-1}) }{   \vert 1-u   \vert} \, d^* u=
\sum_{m=1}^\infty \log p \,\,g(p^m).
\end{equation}
At the archimedean place one has instead
$$
\frac 12\left( \frac{1}{1-\frac 1u}+\frac{1}{1+\frac 1u}\right)=\frac{u^2}{u^2-1}.
$$
The above equation is applied for $u>1$, in which case one can write equivalently
\begin{equation}\label{weil1}
\frac 12\left( \frac{1}{\vert 1-u^{-1}\vert}+\frac{1}{\vert 1+u^{-1}\vert}\right)=\frac{u^2}{u^2-1}.
\end{equation}
Thus,  the term corresponding to \eqref{weil2} yields the distribution $\kappa(u)$ of \eqref{kappadu}.\vspace{.02in}

 In section \ref{hyper}, we shall re-consider the fixed points for the action of the id\`ele class group $C_\K$ on $M=\A_\K/\K^*$. This set is the union of the prime ideals $
 \ffp_{v}=\{ x\in  M \,|\,  x_v=0\}
$, as $v$ varies among the places of the global field $\K$.

\section{The geometry of monoids}\label{sectmonoids}

We denote by $\Mo$  the category of commutative `pointed' monoids $(M,0)$, where $M$ is a semigroup with a commutative multiplicative operation and an identity element $1$. The element $0\in M$ is absorbing for $M$, \ie it satisfies: $0x =x0= 0,~\forall x\in M$.
The morphisms in $\Mo$ are unital homomorphisms of pointed monoids $\varphi: (M,0)\to (N,0)$, thus satisfying the conditions $\varphi(1_M) = 1_N$ and  $\varphi(0) = 0$.\vspace{0.02in}

 The theory of $\Mo$-schemes as in our earlier papers \cite{announc3} and \cite{jamifine} develops, in parallel with the classical  theory of $\Z$-schemes as in \cite{demgab}, the notion of a scheme as  covariant functor from $\Mo$ to the category of sets. This approach is based on the earlier geometric theories over monoids developed  by K. Kato \cite{Kato}, A. Deitmar \cite{deit}, N.~Kurokawa,  H.~Ochiai, M.~Wakayama \cite{KOW},  B.~T\"oen  and M.~Vaqui\'e \cite{TV}.

Several basic  notions holding for (commutative) rings and their geometric counterparts naturally generalize to pointed monoids. From the algebraic side, we recall that an ideal $I$ of a monoid  $M$ is a subset $I\subset M$ such that $0\in I$ and it satisfies the property:
  \begin{equation}\label{ideal}
 x\in I \implies xy\in I\qqq y\in M \,.
\end{equation}
 In particular, an ideal $\ffp\subset M$ is said to be {\em prime} if its complement $\ffp^c=M\setminus \ffp$   is a non-empty multiplicative subset of $M$, \ie if the following condition is satisfied:
 \begin{equation}\label{prime}
  x\notin \ffp, y\notin \ffp \implies xy \notin \ffp.
 \end{equation}
The complement $\ffp_M=(M^\times)^c$ of the invertible elements of $M$ is the largest prime ideal.

From the geometric side, we recall that a geometric monoidal space is a pair $(X,\cO_X)$ where
\vspace{.05in}

$\bullet$~$X$ is a topological space
\vspace{.05in}

$\bullet$~$\cO_X$ is a sheaf of monoids.\vspace{.05in}

The notions of a morphism (of monoidal spaces) as well as those of a prime spectrum $\Sp M$ and  of a geometric scheme are adapted directly from the classical ones (\cf \cite{deit}). The topology on $\Sp M$ admits  a basis of open subsets of the following type  (as $f$ varies in $M$)
 $$
 D(f)=\{\ffp\in \Sp M|\, f\notin \ffp\}.
 $$

By definition, an  $\Mo$-functor is a covariant functor from $\Mo$  to the category of sets. A morphism of $\Mo$-functors, $\phi:X\to Y$, is a natural transformation and as such it determines (a family of) maps of sets
 $$
 \phi_M:X(M)\to Y(M) \qqq M\in\text{obj}(\Mo)
 $$
compatible with any homomorphism $\rho:M\to M'$ in $\Mo$.\vspace{.05in}

A new and interesting property fulfilled by any $\Mo$-functor is that of being automatically {\em local} (\cf\cite{jamifine} \S~3.4.1). Thus, the only requirement that an $\Mo$-functor has to satisfy to define an $\Mo$-scheme is to admit an open covering of affine sub-functors.
We recall that a sub-functor $Y$ of an $\Mo$-functor $X$ (\ie $Y(M)\subset X(M)$ for every object $M$ of $\Mo$) is said to be open in $X$ if for every morphism $\varphi: \Spec(M) \to X$ of $\Mo$-functors, there exists an ideal $I\subset M$ such that for every object  $N$ of $\Mo$ and every $\rho\in\Spec(M)(N)=\Hom_{\Mo}(M,N)$ one has:
$$
\varphi(\rho)\in Y(N)\subset X(N)~\Leftrightarrow~\rho(I)N = N.
$$
Let $X$ be an $\Mo$-functor. A family $\{X_\alpha\}$ of sub-functors of $X$ determines a covering of $X$ if for every abelian group $H$, the following equality (of sets) is fulfilled: $X(\F_1[H])=\bigcup_\alpha X_\alpha(\F_1[H])$.

Any $\Mo$-scheme has an uniquely associated geometric realization (\cf\cite{jamifine} \S~3.4.5). This is a geometric monoidal space whose construction presents several analogies with the geometric realization of a $\Z$-scheme but also a few new interesting features specifically inherent to the discussion with monoids. The most important (new) property states that the full sub-category of the abelian groups inside $\Mo$, whose objects replace the fields within $\Mo$, admits the {\em final object}  $\F_1 = \F_1[\{1\}]$ (\cf\cite{jamifine} \S~3.4.5). This fact implies a remarkable simplification in the description of the geometric realization of an $\Mo$-scheme as stated by the following

\begin{thm}   The geometric space $|X|$  associated to an $\Mo$-scheme $X$ is characterized by  the property
$$
X(M)=\Hom(\Sp\,M, |X|)\qqq M\in\text{obj}(\Mo).
$$
The set underlying $|X|$ is canonically identified with $X(\F_1)$. The topology of $|X|$ is determined by the open subfunctors of $X$ and the structure sheaf by the morphisms to the affine line $\cD$, \ie the functor $\cD(M)=M$ for all  $M\in\text{obj}(\Mo)$.
\end{thm}

We refer to \cite{jamifine} Theorem 3.34 for the detailed statement and proof.

The following canonical projection map  describes a new feature of $\Mo$-schemes which does not hold in general for $\Z$-schemes.  For an $\Mo$-scheme $X$ and for all $M\in\text{obj}(\Mo)$ we define
\begin{equation}\label{projcan}
    \pi_M\,:\,   X(M)\to |X|\,, \ \pi_M(\phi)=\phi(\ffp_M)\qqq \phi\in \Hom_{\Mo}(\Sp(M),|X|)
\end{equation}
where $\ffp_M$ is the largest prime ideal of $M$. Then, for any open subset $U$ of $|X|=X(\F_1)$ with associated sub-functor $\underline U\subset X$  one has
\begin{equation}\label{projcan3}
   \underline  U(M)=\pi_M^{-1}(U)\subset  X(M).
\end{equation}

\subsection{The $\Mo$-scheme $\spadu$}
A basic fundamental  example of an $\Mo$-scheme is the projective line $\spadu$ over $\F_1$. We review shortly its description. As $\Mo$-functor,  $\spadu$ is defined by
  \begin{equation}\label{projspace}
\spadu(M)=M\cup_{M^\times}M,\quad\forall~M\in\text{obj}(\Mo)
\end{equation}
where the gluing map on the disjoint union $M\coprod M\equiv M\times\{1,2\}$ is given by the equivalence relation $(x,1)\sim (x^{-1},2)$, $\forall~x\in M^\times$.

The topological space  $\spad$ of its geometric realization  (\cf \cite{deit}, \cite{announc3})  has three points
\begin{equation}\label{projspace1}
\spad=\{0,u,\infty\}\, , \ \ \overline{\{0\}}=\{0\}\, , \ \overline{\{u\}}=\spad\, , \ \ \overline{\{\infty\}}=\{\infty\}.
\end{equation}
There are three open sets $U_\pm$ and $U=U_+\cap U_-$ in $\spad$
\begin{equation}\label{projspace2}
U_+=\spad\backslash\{\infty\}\,, \ \ U_-=\spad\backslash\{ 0\}\,, \ \ U=U_+\cap U_-.
\end{equation}

\subsection{The monoid $M=\A_\K/\K^\times$ of   ad\`ele  classes}

Let $\K$ be a global field. The product in the ring of ad\`eles over $\K$, descends to the quotient $\A_\K/\K^\times$ to define a natural structure of (commutative) monoid
\begin{equation}\label{monoideclass}
 M=\A_\K/\K^\times, \quad  \K^\times=\GL_1(\K).
\end{equation}
The group $C_\K$ of id\`ele classes is interpreted as the group $M^\times$ of the invertible elements of the monoid $M$. The canonical projection \eqref{projcan} for the $\Mo$-scheme $\spadu$ determines in particular the map ($M = \A_\K/\K^\times$)
\begin{equation}\label{projback}
    \pi_M\,:\,   \spadu(M)=M\cup_{M^\times}M\to \spad.
\end{equation}
$\pi_M$ associates the point $u\in \spad$ to each element of  $M^\times=C_\K$ and either $0$ or $\infty$ to the other elements of  $\spadu(M)$.  We define  the {\em projective ad\`ele class space}  to be the set $\spadu(M)= M\cup_{M^\times}M$.

\subsection{The sheaf $\Omm$  of half-densities on $\spadu(M)$}

To define a natural space $\cS(M)$ of functions on the quotient space $M=\A_\K/\K^\times$  one considers the space of  coinvariants for the action of $\K^\times$ on the Bruhat-Schwartz space $\cS(\A_\K)$ of complex-valued functions on the ad\`eles of $\K$. This action is described explicitly by $f_q(x)=f(qx),~\forall x\in\A_\K,~ q\in\K^\times$. More precisely, one starts with the exact sequence
\begin{equation}\label{fonction1}
0\to \cS(\A_\K)_0 \to \cS(\A_\K)\stackrel{\epsilon}{\to} \C\oplus \C[1]\to 0
\end{equation}
associated to the kernel of the $\K^\times$-invariant linear map $\epsilon(f)=(f(0),\int_{\A_\K} f(x)dx)\in \C\oplus \C[1]$ and then one sets
\begin{equation}\label{fonction3}
\cS(M) := \cS_0(M)\oplus \C\oplus \C[1]\,, \ \  \cS_0(M):=\cS(\A_\K)_0 /\overline{\{f-f_q\}}
\end{equation}
where $\overline{\{f-f_q\}}$ denotes the closure of the sub-space of $ \cS(\A_\K)_0$ generated by the differences $f-f_q$, with $q\in\K^\times$.\vspace{.02in}

We now introduce the functions on the projective ad\`ele class space $\spadu(M)$. The following space of sections determines uniquely a sheaf $\Omm$ on $\spad$, the restriction maps are defined in \eqref{newrestmaps} here below
\begin{eqnarray}\label{sheaf1}
   \Gamma(U_+,\Omm) &=& \cS(M)\nonumber \\
   \Gamma(U_-,\Omm) &=& \cS(M) \nonumber \\
   \Gamma(U_+\cap U_-,\Omm) &=& \strong(C_\K) \nonumber
\end{eqnarray}
where $ \cS(C_\K)$ is the  Bruhat-Schwartz space over $\C_\K$. For a number field $\K$,  $\strong(C_\K)$ is defined as follows
\begin{equation}\label{SCK}
\strong(C_\K)=\,\bigcap_{\beta \in \R}\,\mu^\beta \cS(C_\K)=\{f\in \cS(C_\K)\,|\,\mu^\beta(f)\in \cS(C_\K)\qqq \beta \in \R\}.
\end{equation}
Here, $\mu \in C(C_\K)$ denotes the module  $\mu:C_\K\to\R_+^*$, $\mu^\beta(g)=\mu(g)^\beta$. When $\K$ is a global field of positive characteristic, $\strong(C_\K)$ is the space of Schwartz functions on $C_\K$ with compact support (\cf \cite{CMbook} Definition 4.107 and \cite{Meyer}).
The natural restriction maps $\Gamma(U_\pm, \Omega)\to \Gamma(U_+\cap U_-,\Omega)$  vanish on the components $\C\oplus \C[1]$ of $\cS(M)$, while on $\cS_0(M)$ they are defined as follows:
\begin{eqnarray}\label{newrestmaps}
  (\Res\, f)(g) &=&  \sum_{q\in \K^\times}f(q g) \qqq f\in \cS_0(M)\subset\Gamma(U_+,\Omm)\nonumber    \\
  (\Res\, h)(g) &=& |g|^{-1}\sum_{q\in \K^\times}h(q g^{-1})
  \qqq h\in \cS_0(M)\subset\Gamma(U_-,\Omm).
\end{eqnarray}

\subsection{Spectral realization on $H^1(\spad,\Omm)$}

The following formulas define an action of $C_\K$ on the sheaf $\Omm$. For $\lambda\in C_\K$, define
\begin{eqnarray}\label{rep1}
  \rep_+(\lambda)f(x)  &=& f(\lambda^{-1}x)\qqq f\in \Gamma(U_+,\Omm)\nonumber\\
  \rep_-(\lambda)f(x)  &=& |\lambda|f(\lambda x)\qqq f\in \Gamma(U_-,\Omm) \nonumber \\
   \rep(\lambda)f(x)  &=& f(\lambda^{-1}x)\qqq f\in \Gamma(U_+\cap U_-,\Omm).  \nonumber
\end{eqnarray}
The generator
$
w=
\left(
  \begin{array}{cc}
    0 & 1 \\
    1 & 0 \\
  \end{array}
\right)
$
of the Weyl group $W\subset{\rm PGL}_2$ acts on $C_\K$ by the automorphism $g\mapsto g^w=g^{-1}$, $\forall g\in C_\K$ and this action defines the semi-direct product $N=C_\K\rtimes W$.  Moreover, $w$ acts on $\spad$ by exchanging $0$ and $\infty$.
The action of $w$ on the sheaf $\Omm$ is given by:
\begin{eqnarray}\label{liftw}
  w^{\#} f  &=& f\in \Gamma(U_-,\Omm)\qqq f\in \Gamma(U_+,\Omm)\nonumber\\
 w^{\#} f  &=& f\in \Gamma(U_+,\Omm)\qqq f\in \Gamma(U_-,\Omm) \nonumber \\
   w^{\#} f(g)  &=& |g|^{-1} f(g^{-1})\qqq f\in \Gamma(U_+\cap U_-,\Omm).
\end{eqnarray}
This action defines the  morphism of sheaves $w_\#:\Omm\to w_*\Omm$. The following result is proven in \cite{announc3} (\cf Proposition 5.4)

\begin{prop}\label{thmequnfunct}
There exists a unique action of $N=C_\K\rtimes W$ on the sheaf  $\Omm$ which agrees with \eqref{liftw} on $W$ and restricts on $C_\K$ to the (twist) representation  $\rep[-\frac 12]:=\rep\otimes\mu^{-1/2}$,  where $\mu^{-1/2}$ is seen as representation of $C_\K$.
\end{prop}
\vspace{.05in}

The \v{C}ech complex of the covering $\cU=\{U_\pm\}$ of $\spad$ has two terms
\begin{eqnarray}
  C^0 &=&  \Gamma(U_+,\Omm)\times  \Gamma(U_-,\Omm) \nonumber \\
  C^1 &=& \Gamma(U_+\cap U_-,\Omm). \nonumber
\end{eqnarray}
We introduce the following map $\Sigma: \cS_0(M)\to \cS_\infty(C_\K)$, $\Sigma(f)(x)=\sum_{q\in \K^\times}f(qx)
$. Then, the co-boundary  $\partial:C^0\to C^1$ is given by
\begin{equation}\label{bord}
    \partial(f,h)(g)=\Sigma(f)(g)-|g|^{-1}\Sigma(h)(g^{-1}).
\end{equation}
 Let  $ \alpha$ be a non-trivial character of the additive group $\A_\K/\K$. The lattice $\K\subset \A_\K$ coincides with its own dual. The Fourier transform on $\cS(\A_\K)_0$
\begin{equation}\label{fourier}
F(f)(a)=\int f(x)\alpha(ax)dx
\end{equation}
 becomes canonically defined modulo the subspace $\{f-f_q\}$. For a proof of the following statement we refer to \cite{announc3} (\cf Lemma 5.3 and Theorem 5.5)

\begin{thm}\label{lemreasspec}
The kernel of the co-boundary map $\partial:C^0\to C^1$ coincides with the graph of the Fourier transform  \begin{equation}\label{hzero}
    H^0(\spad,\Omm)=\{(f,F(f))\,|\, f\in \cS(\A_\K)_0/\overline{\{f-f_q\}}\}\oplus \C^{\oplus 2}\oplus (\C[1])^{\oplus 2}.
\end{equation}
The representation $\rep[-\frac 12]$ of $C_\K$ on $H^1(\spad,\Omm)$ determines the spectral realization of the zeros of the $L$-functions. This representation is invariant under the symmetry $\chi(g)\mapsto \chi(g^{-1})$ of the group of Gr\"ossencharakters of the global field $\K$.
\end{thm}

\section{Hyperstructures}\label{hypers}

In this section we briefly recall the results of \cite{wagner} showing that the notion of hyperring  introduced by M. Krasner allows one to understand the algebraic structure of the \ad class space  $\H_\K=\A_\K/\K^\times$ of a global field $\K$.

\subsection{Hypergroups and hyperrings}

We start by reviewing the notion of a {\em canonical} hypergroup $(H,+)$. For our applications it will be enough to consider the commutative case and we  denote by $+$ the hyper-composition law in $H$. The novelty is that for hypergroups such as $H$ the sum $x+y$ of two elements in $H$ is no longer a single element of $H$ but a non-empty subset of $H$. It is customary to  define a hyper-operation on $H$ as a map
\[
+: H\times H \to \mathcal P(H)^*
\]
 taking values into the set $\cP(H)^*$ of all non-empty subsets of $H$.  One uses the notation $\forall A,B\subseteq H,~A+B:=\{\cup (a+b)~|a\in A, b\in B\}$. The definition of a commutative canonical hypergroup requires that $H$ has a neutral element $0\in H$ (\ie an additive identity) and that the following axioms apply:\vspace{.05in}

 $(1)$~$x+y=y+x,\qquad\forall x,y\in H$\vspace{.05in}

$(2)$~$(x+y)+z=x+(y+z),\qquad\forall x,y,z\in H$ \vspace{.05in}

$(3)$~$0+x=x= x+0,\qquad \forall x\in H$\vspace{.05in}

$(4)$~$\forall x  \in H~  \ \exists!~y(=-x)\in H\quad {\rm s.t.}\quad 0\in x+y$\vspace{.05in}

$(5)$~$x\in y+z~\Longrightarrow~ z\in x-y.$\vspace{.05in}

 Property $(5)$ is usually called {\em reversibility}.
 \vspace{.05in}

\begin{lem}\label{general}  Let $(G,\cdot)$ be a commutative  group, and let $K\subset \Aut(G)$ be a subgroup of the group of automorphisms of $G$. Then the following operation defines a structure of hypergroup on the set $H=\{K(g)|g\in G\}$ of the orbits of the action of $K$ on $G$:
\begin{equation}\label{autg}
    K(g_1)\cdot K(g_2):=(Kg_1\cdot Kg_2)/ K.
\end{equation}
\end{lem}

The notion of  a  {\em hyperring}  (\cf \cite{Krasner}, \cite{Krasner1}) is the natural generalization of the classical notion of a ring, obtained by replacing a classical additive law by a hyperaddition. More precisely, a hyperring  $(R,+,\cdot)$ is a non-empty set $R$ endowed with a hyperaddition $+$ and a multiplicative operation $\cdot$ satisfying the following properties:\vspace{.05in}

$(a)$~$(R,+)$ is a commutative canonical hypergroup\vspace{.05in}

$(b)$~$(R,\cdot)$ is a monoid with multiplicative identity $1$\vspace{.05in}

$(c)$~$\forall r,s,t\in R$:~~$r(s+t) = rs+rt$ and $(s+t)r = sr+tr$\vspace{.05in}

$(d)$~$\forall r\in R$:~~$r\cdot 0=0\cdot r =0$, \ie $0\in R$ is an absorbing element\vspace{.05in}

$(e)$~$0\neq 1$.\vspace{.05in}

Let $(R_1,+_1,\cdot_1)$, $(R_2,+_2,\cdot_2)$ be two hyperrings. A map $f: R_1 \to R_2$ is called a homomorphism of hyperrings if the following conditions are satisfied\vspace{.05in}

$(1)$~$f(a+_1 b)\subseteq f(a)+_2 f(b),~\forall a,b\in R_1$\vspace{.05in}

$(2)$~$f(a\cdot_1 b) = f(a)\cdot_2 f(b),~\forall a,b\in R_1.$\vspace{.05in}

\subsection{$\kras$, $\sign$ and the \ad class space $\H_\K=\A_\K/\K^\times$}

A hyperring $(R,+,\cdot)$ is called a {\em hyperfield} if $(R\setminus\{0\},\cdot)$ is a group. The most basic example of a hyperfield is the {\em Krasner hyperfield} $\kras=(\{0,1\},+,\cdot)$ with additive neutral element $0$, satisfying the hyper-rule: $1+1=\{0,1\}$ and with the usual   multiplication, with identity $1$. Likewise $\F_2$ encodes the arithmetic of even and odd numbers, $\kras$ encodes the arithmetic of zero and non-zero numbers. The hyperfield $\kras$ is the natural extension, in the category of hyperrings, of the commutative (pointed) monoid $\F_1$, \ie  $(\kras,\cdot) = \F_1$.

Another interesting example of basic hyperstructure is the {\em hyperfield of signs} $\sign=(\{0,\pm 1\},+,\cdot)$ where the hyper-addition is given by the ``rule of signs"
\begin{equation}\label{addsign}
    1+1=1\,, \ -1-1=-1\,, \ 1-1=-1+1=\{-1,0,1\}
\end{equation}
and where the usual multiplication is also given by the rule of multiplication of signs. $\sign$ encodes the arithmetic of the signs of numbers and it is the natural extension, in the category of hyperrings, of the commutative (pointed) monoid $\F_{1^2}$, \ie  $(\sign,\cdot) = \F_{1^2}$.

There is a  unique hyperring homomorphism $\sigma: \Z\to \sign$, $\sigma(n) = \text{sign}(n),~\forall n\neq 0$, $\sigma(0) = 0$. Moreover, the absolute value determines a canonical surjective homomorphism of hyperfields $\pi: \sign \to \kras$. Then, by using the composite homomorphism $h = \pi\circ \sigma$, one can perform the extension of scalars from $\Z$ to $\kras$ and show that for any commutative ring $R$ containing $\Q$ as a subfield, one obtains the isomorphisms $R\otimes_\Z\kras = R/\Q^\times$ and $R\otimes_\Z\sign = R/\Q_+^\times$ (\cf\cite{wagner} Proposition~6.1).\vspace{.02in}

Let $R$ be a commutative ring and let $G\subset R^\times$ be a subgroup of its multiplicative group. Then
  the following operations define a hyperring structure on the set $R/G$ of orbits for the action of $G$ on $R$ by multiplication\vspace{.05in}

  $\bullet$~Hyperaddition
  $$
x+ y:=\left(xG+yG\right)/G \qqq x,y\in R/G
$$

$\bullet$~Multiplication
$$
xG\cdot yG=xyG\qqq x,y\in R/G.
$$

In particular, one may start with a field $K$ and consider the hyperring  $K/K^\times$. This way, one obtains a hyperstructure whose underlying set is made by two classes \ie the class of $0$ and that of $1$. It is easy to see that if $K$ has more than two elements, then $K/K^\times$ is isomorphic to the Krasner hyperfield $\kras$ .

In general, a hyperring need not  contain $\kras$ as a sub-hyperfield. For quotient hyperrings, like the ones we have introduced right above, there is a precise condition that ensures the occurrence of this case (\cf\cite{wagner} Proposition 2.6)
\begin{thm}\label{wagnad} Let $R$ is a commutative ring and $G\subset R^\times$   a proper subgroup of the multiplicative group of units of $R$, then,
  the hyperring $R/G$ contains $\kras$ as a sub-hyperfield if and only if $\{0\}\cup G$ is a subfield of $R$.
  \end{thm}
Since the \ad class space $\H_\K=\A_\K/\K^\times$ of a global field $\K$ is the quotient of the commutative ring $R=\A_\K$ by $G=\K^\times$ and $\{0\}\cup G=\K$ is a subfield of $R=\A_\K$, one obtains
\begin{cor}\label{wagnad1} The \ad class space $\H_\K=\A_\K/\K^\times$ of a global field $\K$ is a hyperring extension of $\kras$.
\end{cor}

It is elementary to prove that for any ring $R$, the map
\[
\varphi:\Spec(R)\to \Hom(R,\kras)\,, \qquad\varphi(\ffp)=\varphi_\ffp
\]
$$
    \varphi_\ffp(x)=0\qqq x\in \ffp\, , \ \ \varphi_\ffp(x)=1\qqq x\notin \ffp
$$
determines a  natural bijection of sets. This fact shows that the hyperfield $\kras$ plays, among hyperrings, the equivalent role of the monoid $\F_1$ among monoids (\cf \cite{jamifine} Proposition 3.32).

\subsection{Extensions of $\kras$ and incidence groups}

In this section we  outline several results which show that the structure of hyperfield and hyperring extensions of $\kras$ is intimately connected to the geometric notion of incidence group of Ellers and Karzel (\cf\cite{Ellers}).    We refer to  \S 3 of \cite{wagner} to read more details and for the proofs.\vspace{.02in}

There is a canonical  correspondence, originally established by  Prenowitz in \cite{Prenowitz},
 between $\kras$-vector spaces $E$ and projective geometries $(\cP,\cL)$ in which every line has {\em at least $4$ points}.
The line passing through two distinct points $x,y$ of $\cP:=E\setminus \{0\}$ is defined by
  $$
    L(x,y)=(x+y)\cup \{x,y\}.
  $$

Conversely, the hyper-addition in $E:=\cP\cup\{0\}$ is defined by the rule
$$x+y = L(x,y)\setminus\{x,y\},\quad\text{if}~~x\neq y,\quad x+x = \{0,x\}.$$

If a group $G$ is the set of points of a projective geometry, then $G$ is called a {\em two-sided incidence group} if the left and the right translations by $G$ are automorphisms of the geometry.

 Let $\H\supset \kras$ be a  hyperfield extension of $\kras$ and let $(\cP,\cL)$ be the associated geometry. Then, the multiplicative group $\H^\times$, endowed with the geometry $(\cP,\cL)$, is a two-sided incidence group.
Conversely, let $G$ be a two-sided incidence group. Then, there exists a unique hyperfield extension $\H\supset \kras$ such that $\H=G\cup \{0\}$.

The  classification of Desarguesian commutative incidence groups due to H. Karzel (\cf\cite{Karzel}) applies to commutative hyperfield extensions $H$ of $\kras$ such that $\dim_\kras\H>3$.
 Let $\H\supset \kras$ be a commutative hyperfield extension of $\kras$. Assume that the geometry associated to the $\kras$-vector space $\H$ is Desarguesian (this condition is automatic  if  $\dim_\kras\H>3$)
and of dimension at least $2$. Then, there exists a  {\em unique} pair $(F,K)$ of a commutative field $F$ and a subfield $K\subset F$ such that
$$
    \H=F/K^\times.
$$

In view of the result just stated, the classification of all {\em finite, commutative} hyperfield extensions of $\kras$ reduces to the determination of non-Desarguesian, finite projective planes with a simply transitive abelian group of collineations. More precisely, if
 $\H\supset \kras$ is a finite commutative  hyperfield extension of $\kras$, then,  one of the following cases occurs:\vspace{.05in}

$(1)$~$\H=\kras[G]$,  for a finite abelian group $G$.\vspace{.05in}

$(2)$~There exists a finite field extension $\F_q\subset \F_{q^m}$ of a finite field $\F_q$ such that $\H=\F_{q^m}/\F_q^\times$.\vspace{.05in}

$(3)$~There exists a finite, non-Desarguesian projective plane $\cP$ and a simply transitive abelian group $G$ of collineations of $\cP$, such that $G$ is the commutative incidence group associated to $\H$.

There are no known examples of finite, commutative hyperfield extensions $\H\supset\kras$ producing  projective planes as in case (3). In fact, there is a conjecture  based on some results of  A. Wagner (\cf\cite{Beutel}, \cite{Wagner}, \cite{Wagner1})  stating that such case cannot occur.  M. Hall proved two results (\cf\cite{Hall}) which imply the following conclusions:\vspace{.05in}

$\bullet$~Assume that $\H^\times$ is cyclic. Let $n+1$ be the cardinality of each line of the geometry. Then for each prime divisor of $n$,  the map $x\mapsto x^p$ is an automorphism of $\H\supset \kras$.\vspace{.05in}

$\bullet$~There exists an {\em infinite} hyperfield extension $\H\supset \kras$ whose geometry is non-Desarguesian and $\H^\times \simeq \Z$.
 \vspace{.05in}

Let $\H\supset \kras$ be a commutative hyperring extension of $\kras$. Assume that $\H$ has no zero divisors and that $\dim_\kras\H >3$. Then, there exists a unique pair $(A,K)$ of a commutative integral domain $A$ and a subfield $K\subset A$ such that
$$
    \H=A/K^\times.
$$

Let $A_j$ ($j=1,2$) be  commutative algebras over two fields $K_j\neq \F_2$ and let
$$
\rho\,:\, A_1/K_1^\times\to A_2/K_2^\times
$$
be a homomorphism of hyperrings. Assume that the range of $\rho$ is of $\kras$-dimension at least $3$, then $\rho$ is induced by a unique ring homomorphism $\tilde \rho:A_1\to A_2$ such that $\alpha=\tilde \rho|_{K_1}$ is a field inclusion $\alpha:K_1\to K_2$. These results show that, in higher rank, the category of hyperring extensions of $\kras$ is the category of algebras over fields
with twisted morphisms.

\section{The hyperfield $\rco$}\label{rconvex}

In this section  we prove that the set of the  real numbers is endowed with a natural structure of hyperfield extension $\R^{\rm convex}$  of the hyperfield of signs $\sign$. It turns out that the hyperstructure on $\R^{\rm convex}$ is a refinement of the algebraic structure on the semi-field $\rmax$ commonly used in idempotent analysis and tropical geometry. The hyperfield $\R^{\rm convex}$  has characteristic one and it comes equipped with a one parameter group of automorphisms which plays the role of the  Frobenius in characteristic one.\vspace{.05in}

\subsection{Sign-convex subsets of $\R$}\label{sectsconvex}

The sign of a real number determines a canonical surjective map
\begin{equation}\label{signmap}
    \ss: \R\to \{0,\pm 1\},\quad \ss(r) = \left\{
                                     \begin{array}{ll}
                                       0, & \hbox{if $r=0$;} \\
                                       1, & \hbox{for $r>0$;} \\
                                       -1, & \hbox{for $r<0$.}
                                     \end{array}
                                   \right.
\end{equation}

For any pair of real numbers $x,y\in\R$, we set
\begin{equation}\label{rconv}
c(x,y)=    \{\alpha x+\beta y\mid \alpha > 0,\beta > 0,\ \ss(\alpha x+\beta y)=\alpha\, \ss(x)+\beta \,   \ss(y)\}.
\end{equation}
\begin{defn} A subset $C\subset \R$ is said to be sign-convex if $\forall x,y\in C$ one has
$c(x,y)\subset C$.
\end{defn}

We {\em fix}  a homeomorphism $\phi: (0,\infty)\to\{e^{i\theta}\mid \theta \in (0,\pi)\}$ of the positive real line with the upper-half unit circle in $\C$, such that $\lim_{x\to 0}\phi(x)= -1$  and
$\lim_{x\to\infty}\phi(x)= 1$.

We let $U=\{0\}\cup\{z\in \C\mid |z|=1,\ z\notin \R\}$
and we extend uniquely $\phi$ to a bijection $\phi: \R\stackrel{\sim}{\to} U$ by setting:
\begin{equation}\label{bij}
    \phi(x)=-\phi(-x)\qqq x\neq 0\,, \ \phi(0)=0.
\end{equation}
For instance one can take
$$
\phi(x)=\ss(x)e^{i\theta(x)}\,, \ \theta(x)=\frac{\pi}{1+x^2}\,.
$$

   \begin{figure}
\begin{center}
\includegraphics[scale=0.7]{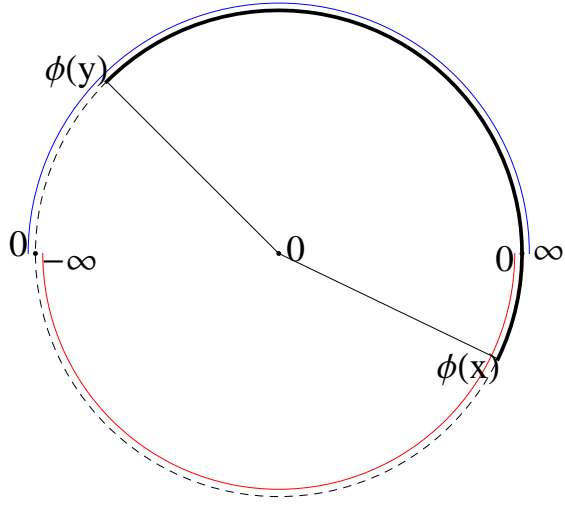}
\end{center}
\caption{The set $\phi(c(x,y))$ for $x<0$, $y>-x$.\label{rconvdraw} }
\end{figure}

\begin{lem}\label{theeq} For all $x,y\in \R$ one has
\begin{equation}\label{equiv1bis}
    z\in c(x,y)\Leftrightarrow \phi(z)\in \R_+^\times\phi(x)+\R_+^\times\phi(y).
\end{equation}
\end{lem}
\proof For $x=0$, one has $c(0,y)=\{y\}$, $\forall y\in\R$ and \eqref{equiv1bis} holds. If $y=0$, one concludes in a similar way. Thus we can assume $x\neq 0\neq y$. Next, we show that
\begin{equation}\label{simp}
0\in c(x,y) \Leftrightarrow y=-x.
\end{equation}
If $0\in c(x,y)$, then there exist $\alpha>0,\beta>0$ such that $\alpha x+\beta y=0$ and
$\alpha\, \ss(x)+\beta \,   \ss(y)=0$. This implies $\alpha=\beta$ and  $y=-x$. Conversely, if $y=-x$ one has $0=x+y$ and $0=\ss(x)+\ss(y)$ so that $0\in c(x,y)$.

Since $0\in \R_+^\times\phi(x)+\R_+^\times\phi(y)$ if and only if $y=-x$, this proves that \eqref{equiv1bis} holds for $z=0$.
We can then assume that $x,y,z$ are all different from $0$. If $0<x<y$  both sides of \eqref{equiv1bis} give the interval $(x,y)$ and one gets \eqref{equiv1bis}  when $x$ and $y$ have the same sign. For $x<0$ and $y>0$ there are three possible cases and using \eqref{rconv} one gets\vspace{.05in}

$\bullet$~If $-x<y$ then $c(x,y)=(x,0)\cup (y,\infty)$.\vspace{.05in}

$\bullet$~If $-x=y$ then $c(x,y)=\{x,0,y\}$.\vspace{.05in}

$\bullet$~If $-x>y$ then $c(x,y)=(-\infty,x)\cup (0,y)$.\vspace{.05in}

Thus \eqref{equiv1bis} can be checked directly in each case: \cf Figure \ref{rconvdraw}.
  \endproof

A subset  $\Gamma\subset \C$ is called a {\em convex cone} if $\Gamma$ is  stable both for addition and for the action of $\R_+^\times$ on $\Gamma$ by multiplication. For any subset $X\subset \C$, the convex cone $\Gamma(X)$ generated by $X$ verifies the equality
\begin{equation}\label{cvx}
    \Gamma(X)\backslash \{0\}=\{\alpha x+\beta y\,, \ \alpha>0, \ \beta>0,\ x,y\in X\}\backslash \{0\}.
\end{equation}
It may happen that $0\in \Gamma(X)$ cannot be written as a sum of two elements of $\R_+^\times X$, but it is always possible to write $0$ as a sum of three elements of $\R_+^\times X$.\vspace{.02in}

Next results shows that sign-convex subsets of $\R$ are determined by convex-cones in $\C$.

\begin{cor}\label{equiv1} Let $\phi$ be as in \eqref{bij} and  $C\subset \R$. The following conditions are equivalent: \vspace{.05in}

$\bullet$~$C$ is sign-convex.\vspace{.05in}

$\bullet$~$C=\phi^{-1}(\Gamma)$ where $\Gamma$ be the cone generated by $\phi(C)$.\vspace{.05in}

$\bullet$~ There exists a convex cone $\Gamma\subset \C$ such that $C=\phi^{-1}(\Gamma)$. \end{cor}
\proof   Let $C\subset \R$ be sign-convex. Let $\Gamma=\Gamma(\phi(C))$ be the cone generated by $\phi(C)$. By \eqref{cvx} any non-zero element of $\Gamma$ is of the form
\begin{equation}\label{cone}
    \xi=\alpha\phi(x)+\beta\phi(y),\quad \text{for some}~\alpha>0, \ \beta>0,\ x,y\in C.
\end{equation}
Thus, by applying  Lemma \ref{theeq}, one has $0\neq\xi=\phi(z)\in\Gamma \implies z\in C$, for $z\neq 0$.
If $0\in \Gamma$, then there exist three elements $x,y,z$ of $C$ and three positive real numbers $\alpha,\beta,\gamma>0$ such that $\alpha\phi(x)+\beta \phi(y)+\gamma \phi(z)=0$.
Since $C$ is sign convex, Lemma \ref{theeq} shows that there exists $a\in C$, with $-a\in C$ so that $0\in c(a,-a)\subset C$.

If $C=\phi^{-1}(\Gamma)$ for some convex cone $\Gamma\subset \C$, then, it follows again from Lemma \ref{theeq} that $C$ is sign-convex. \endproof

If $C\subset \R$ is sign-convex then the two subsets $C^\pm=C\cap \pm(0,\infty)\subset\R$ are convex \ie they are intervals, but the converse of this statement fails.  If $C$ is sign-convex, so is $C\cup\{0\}$. Moreover, set aside $\R$, the only sign-convex subsets $C\subset\R$ which contain a pair $x,-x$,  for some $x>0$ are of the following types\vspace{.05in}

$\bullet$~$C_x=\{-x,0,x\}$\vspace{.05in}

$\bullet$~$C_x^+=[-x,0]\cup[x,\infty)$\vspace{.05in}

$\bullet$~$C_x^-=(-\infty,-x]\cup [0,x]$.\vspace{.05in}

\begin{lem}\label{signhull} Let $C\subset \R$ be sign-convex and $x\notin C$. Then the smallest sign-convex set containing $x$ and $C$ is the set
$$
C'=\{x\}\cup C\cup_{y\in C}c(x,y).
$$
\end{lem}
\proof
Since any sign-convex subset containing $x$ and $C$ contains $C'$, and since $C'$ contains $x$ and $C$, it suffices to show that $C'$ is sign-convex. This follows from Corollary \ref{equiv1}. \endproof

\begin{lem}\label{theeqbis} For all $x,y,z\in \R$ one has
\begin{equation}\label{equiv1ter}
    t\in \bigcup_{u\in c(x,y)}c(u,z)\Leftrightarrow \phi(t)\in \R_+^\times\phi(x)+\R_+^\times\phi(y)+\R_+^\times\phi(z).
\end{equation}
\end{lem}
\proof By Lemma \ref{theeq} one has $\phi(t)\in \R_+^\times( \R_+^\times\phi(x)+\R_+^\times\phi(y))+\R_+^\times\phi(z)$ for any $t\in \cup_{u\in c(x,y)}c(u,z)$. Conversely, let $t\in \R$ be such that
$$
\phi(t)=\lambda_1\phi(x)+\lambda_2\phi(y)+\lambda_3\phi(z)\,, \ \lambda_i>0\,.
$$
Let $\alpha=\lambda_1\phi(x)+\lambda_2\phi(y)$.
If there exists $\lambda>0$ such that $\lambda \alpha \in U$ (\cf\eqref{bij}) then there exists $u\in \R$
with $\phi(u)=\lambda \alpha\in \R_+^\times\phi(x)+\R_+^\times\phi(y)$. Thus by Lemma \ref{theeq} one has $u\in c(x,y)$ and since
$$
\phi(t)=\lambda^{-1}\phi(u)+\lambda_3\phi(z)\in \R_+^\times\phi(u)+\R_+^\times\phi(z)
$$
one gets $t\in c(u,z)\subset \cup_{v\in c(x,y)}c(v,z)$ as required. Otherwise one has $\alpha\in \R$, $\alpha\neq 0$.  One has $\phi(x)\notin \{0,\pm \phi(y)\}$. Let $\epsilon=\ss(\alpha)\in\{-1,1\}$. Then there exists an open neighborhood $V$ of $\epsilon$ in the unit circle such that
\begin{equation}\label{Vsubset}
V\subset \R_+^\times\phi(x)+\R_+^\times\phi(y).
\end{equation}
One has $\phi(t)=\alpha+\lambda_3\phi(z)$ by construction, and since $\phi(z)\neq \epsilon$ it follows that $\phi(t)$ is in the interior of the short interval  between
$\epsilon$ and $\phi(z)$ on the unit circle. Thus there exists $v\in V$, $v\neq \epsilon$ such that $\phi(t)$ is in the interior of the short interval  between
$v$ and $\phi(z)$ on the unit circle. Let $u\in \R$ such that $v=\phi(u)$, then by \eqref{Vsubset} and Lemma \ref{theeq}, one has $u\in c(x,y)$. Moreover since
$\phi(t)$ is in the interior of the short interval  between
$\phi(u)$ and $\phi(z)$ one has $t\in c(u,z)$ again by Lemma \ref{theeq}.\endproof

\subsection{Construction of $\R^{\text{convex}}$}

We now use the preliminary results of the previous subsection to construct the hyperfield extension $\rco$ of $\sign$.

\begin{thm} \label{rconvprop}On the set $\R$ there exists a unique structure of hyperfield $\rco=(\R, +_c,\cdot)$, where  $\forall x,y\in \R$ one sets $x+_c y=c(x,y)$ \ie
\begin{equation}\label{rcosum}
    x+_c y:=\{\alpha x+\beta y\mid \alpha>0,\beta>0,\ \ss(\alpha x+\beta y)=\alpha\, \ss(x)+\beta \,   \ss(y)\}
\end{equation}
and where the multiplication $\cdot$ is the classical one. The hyperaddition $+_c$ on $\rco$ is uniquely determined by the properties:\vspace{.05in}

$(1)$~$x+_c y = (x,y)$\quad$\forall~y>x>0$\vspace{.05in}

$(2)$~$\sign \subset\R^{\text{convex}}$ as a sub-hyperfield.

\end{thm}

\proof The operation $x+_c y=c(x,y)$ is commutative by construction. For $x=0$ one has
$$
0+_cy=\{\beta y\mid \beta>0,\ss(\beta y)=\beta \,   \ss(y)\}=\{y\}.
$$
Thus $0$ is a neutral element. By \eqref{simp}, one has $0\in x+_c y \Leftrightarrow y=-x$.
 Note that one has
\begin{equation}\label{xplus}
    x+_c (-x)=\{-x,0,x\}\qqq x\in \R
\end{equation}
and also
\begin{equation}\label{xplusx}
    x+_c x=\{x\}\qqq x\in \R.
\end{equation}
Moreover, we claim that for any real number $a$ one has
\begin{equation}\label{distr}
    a(x+_c y)=ax+_c ay.
\end{equation}
This holds for $a=0$. For $a>0$ the statement follows from \eqref{rcosum}, by using the equality $\ss(az)=\ss(z)$ which is valid for all $z\in \R$. The claim holds also for $a<0$, in fact it follows by applying \eqref{rcosum} and $\ss(az)=-\ss(z)$ for all $z\in \R$.

The associativity of the hyper-addition $+_c$ follows from Lemma \ref{theeqbis} which shows that $(x+_cy)+_cz=\cup_{u\in c(x,y)}c(u,z)$ is symmetric in $x,y,z$.
Now, we prove the reversibility  of the hyper-addition $+_c$. By Lemma \ref{theeq} one has
$$
-z\in x+_c y \Leftrightarrow 0\in \R_+^\times\phi(x)+\R_+^\times\phi(y)+\R_+^\times\phi(z).
$$
Since this equivalent condition is symmetric in $x,y,z$ we obtain the reversibility.

To prove the uniqueness of $\rco$ it is enough to determine the  hypersum $1+_c x\subset\R$ assuming $x\notin\{0,\pm 1\}$. For $x>0$, this is the interval between $1$ and $x$. Assume then $x<0$. For $y<0$, one has using reversibility,
\[
y\in 1+_c x~\Leftrightarrow~x\in y-_c 1~\Leftrightarrow~-x\in -y+_c 1 = \{-\lambda y+(1-\lambda)|\lambda\in(0,1)\}.
\]
Thus, for $x<0$, we know the description the intersection $(1+_c x)\cap(-\infty,0)$, \ie if  $|x|<1$, it is the interval $(x,0)$ and if $|x|>1$ it is $(-\infty,x)$. By applying the distributivity  we know, for $x<0$, that $(-1-_cx)\cap(0,\infty)=-\left((1+_cx)\cap(-\infty,0)\right)$, and hence this determines $t(-1-_cx)\cap(0,\infty)$ for $t>0$. Taking $t=-1/x>0$, this determines $(1/x+_c1)\cap(0,\infty)$ which gives
\begin{align*}
\text{if $|x|<1$}~&:~t(-(x,0)) = t(0,-x)=(0,1)\\
\text{if $|x|>1$}~&:~t(-(-\infty,x)) = t(-x,\infty)=(1,\infty).
\end{align*}
Thus, for $x<0$ we get, replacing $x$ by $1/x$,
\begin{align*}\label{16j}
|x|<1,~x<0~&\Rightarrow~1+_c x=(x,0)\cup(1,\infty)\\
|x|>1,~x<0~&\Rightarrow~1+_c x=(-\infty,x)\cup(0,1)
\end{align*}
which gives the required uniqueness.\endproof

\begin{prop} \label{autprop} Let $\Aut(\rco)$  be the group of automorphisms of the hyperfield $\rco$. The following map defines an isomorphism of (multiplicative) groups $\R^\times\stackrel{\sim}{\to} \Aut(\rco)$,  $\lambda\mapsto \theta_\lambda$, where
\begin{equation}\label{autr}
    \theta_\lambda(x) = x^\lambda\quad \forall x >0, \quad  \theta_\lambda(x) = x\qqq x\in \sign \subset \rco.
\end{equation}
\end{prop}
\proof Let first check that \eqref{autr} defines an automorphism $\theta_\lambda\in \Aut(\rco)$. Since $\theta_\lambda$ is compatible with the product, one just needs to show the compatibility with the hyperaddition $+_c$. This can be checked directly.
Let then $\alpha \in \Aut(\rco)$ be an automorphism. Since $-1$ is the unique additive inverse of $1$, one gets that $\alpha(x) = x\qqq x\in \sign \subset \rco$. Since the subgroup $\R_+^\times$ of the multiplicative group is the subgroup of squares, it is preserved globally by $\alpha$ and thus $\alpha$ defines a group automorphism of $\R_+^\times$. Furthermore, since $\alpha$ is compatible with the hyperaddition its restriction to $\R_+^\times$ is monotonic and  hence it is given by $\alpha(x)=x^\lambda$ for some $\lambda \in \R^\times$.
\endproof

\begin{rem}{\rm There is no intermediate sub-hyperfield $\sign\subset F\subset \rco$, since for $\xi\notin\{0,\pm 1\}$, one would have $\pm\xi>0$ and then $1+\xi$ would contain an open interval generating the multiplicative group (one has $-1\in F$).}
\end{rem}

\subsection{The hyperfields $\sign[G]$}

In this section we show how to extend the construction of $\rco$ to a functor
 $G\mapsto \sign[G]$ from the category of totally ordered non-discrete abelian groups to hyperfield extensions of $\sign$.
For $\kras$-vector spaces $E$ with $\dim_\kras(E)=2$, Remark~3.7 of \cite{wagner} shows a simple  ``set-theoretic'' construction of this functor since in that case there is a single line in the corresponding geometry. We first describe an analogous construction for a $\sign$-vector space, starting with an ordered set $G$ such that ($a,b,c\in G$)
\begin{equation}\label{7j}
\forall a<b,~\exists c,~a<c<b,\quad G\ \text{has no minimal or maximal element.}
\end{equation}
On the set $$\sign(G) := -G\coprod\{0\}\coprod G$$  we define the hyperaddition as follows:
\begin{xalignat}{2}\label{8j}
x&>0,~y>0,~~x\neq y &\quad x+y &:= \{z|\inf(x,y)<z<\sup(x,y)\}\\\notag
x&=y>0 &\quad  x+x &:=x\\\notag
x&>0 &\quad  x+(-x)&:=\{0,\pm x\}\\\notag
x&>0,~y<0,~x>-y &\quad  x+y &:= \{z|z>x\}\cup-\{t| t<|y|\}.\notag
\end{xalignat}
Then, we extend this hyperaddition uniquely using the rule $-(a+b)=-a-b,~~0+x=x+0 = x$.

\begin{prop} $(\sign(G),+)$ is a $\sign$-vector space.
\end{prop}

The main point in the proof of the above proposition is to check the associativity and then to note that all ordered configurations which occur in $\sign(G)$ also occur in  $\rco$. Thus the associativity follows from that case using the condition \eqref{7j} to obtain the symmetry of $(x+y)+z$  in $x,y,z$.

\begin{thm} Let $G$ be a totally ordered commutative group fulfilling \eqref{7j}, then $\sign(G)$ is a hyperfield extension of $\sign$.
\end{thm}
\proof It follows from the above proposition that $\sign(G)$ is a hypergroup. Moreover the operation $x\to -x$ is an automorphism of $\sign(G)$. We denote $G$ as a multiplicative group. The multiplication $L_a$ by a fixed element  $a\in G$ preserves the order and hence it is compatible with the hyperaddition $L_a(x+y) := L_a(x)+L_a(y)$. This shows that if one endows $\sign(G)$ with the product of $G$ satisfying $(-1)^2 = 1$, then we get a hyperfield. One has $1+1=1$ and $-1+1 = \{0,\pm 1\}$, thus $\sign(G)$ contains the hyperfield of signs $\sign$.
\endproof
By applying this result to $G = \R_+^\times$ ($\sim \R$) as ordered multiplicative group, one then gets a canonical isomorphism
\begin{equation}\label{14j}
\sign(\R_+^\times) \simeq \rco\,.
\end{equation}

\section{Function theory for hyperrings}\label{fnthhr}

We recall that   a function on a scheme $X$, viewed as a covariant $\Z$-functor $\underline X:\An\to\Se$, is  a morphism of $\Z$-functors $f:\underline X\to \cD$, where $\cD$ is the functor affine line $\cD=\mathfrak{spec}(\Z[T])$,  with geometric scheme  $\Spec(\Z[T])$ (\cf \cite{demgab} Chapter I and \cite{announc3}). When $X=\Spec(R)$, with $R\in\text{obj}(\An)$ (\ie $R$ a commutative ring with unit), one derives a natural identification of functions on $X$ with elements of the ring $R$
\begin{equation}\label{ztidd}
   \Hom_\An(\Z[T],R)\simeq R.
\end{equation}
 In the category of hyperrings, the identification \eqref{ztidd} no longer holds in general. Indeed, $\kras$ has only two elements while $\Hom_\han(\Z[T],\kras)\simeq \Spec(\Z[T])$ is countably infinite (\cf\cite{wagner} Proposition~2.13). In the following sections we take up the study of function theory on the spaces $\Spec \kras$ and  $\Spec \sign$ \ie we describe the sets $\cD(\kras)=\Hom_\han(\Z[T],\kras)$ and of $\cD(\sign)=\Hom_\han(\Z[T],\sign)$
 together with the hyperoperations coming from addition and multiplication of functions.

 \subsection{Coproducts and homomorphisms to hyperrings}\label{coprosect}

  Let $\cH$ be a commutative ring with unit, and let $\Delta:\cH\to \cH\otimes_\Z\cH$ be a coproduct. Given two ring homomorphisms $\varphi_j: \cH\to R$ ($j=1,2$) to a commutative ring $R$, the composition
$\varphi=(\varphi_1\otimes \varphi_2)\circ \Delta$ defines a ring homomorphism $\varphi: \cH\to R$. When $R$ is a hyperring,
one introduces the following notion
\begin{defn} Let $(\cH,\Delta)$ be a commutative ring with a coproduct and let $R$ be a hyperring. Let $\varphi_j\in \Hom_\han(\cH,R)$, $j=1,2$. One defines
\begin{multline}\label{equhopf}
\varphi_1\star_\Delta \varphi_2 = \{\varphi\in \Hom_\han(\cH,R)|\varphi(x) \in \sum \varphi_1(x_{(1)})\varphi_2(x_{(2)}) ,\\ \text{for all decompositions}\ \ \Delta(x)=\sum x_{(1)}\otimes x_{(2)} \}.
\end{multline}
\end{defn}
In general, for $x\in\cH$, there are several ways  to write
\begin{equation}\label{5h}
\Delta(x)=\sum x_{(1)}\otimes x_{(2)}
\end{equation}
which  represent  the {\em same} element of $\cH\otimes\cH$. The condition
$\varphi(x) \in \sum \varphi_1(x_{(1)})\varphi_2(x_{(2)})$ has to hold for all these decompositions.
 In general, $\varphi_1\star_\Delta \varphi_2$ can be empty or it may contain several elements. When $\varphi_1\star_\Delta \varphi_2=\{\varphi\}$ is made by a single element we simply write  $\varphi_1\star_\Delta \varphi_2= \varphi $.\vspace{.05in}

 The canonical homomorphism $\Z\to \cH$ induces a restriction homomorphism
 \begin{equation}\label{pi}
   \pi:\Hom(\cH,R)\to \Hom(\Z,R).
 \end{equation}
 The following lemma shows that $\varphi_1\star_\Delta \varphi_2$ is empty when the restrictions $\pi(\varphi_j)=\varphi_j|_\Z$ are distinct.

  \begin{lem} \label{restoz} Let $(\cH,\Delta)$ be a commutative ring with a coproduct and let $R$ be a hyperring.  Let $\varphi_j\in \Hom_\han(\cH,R)$, $j=1,2$. If $\varphi\in\varphi_1\star_\Delta\varphi_2$, then
\begin{equation}\label{14h}
\varphi|_\Z=\varphi_1|_\Z=\varphi_2|_\Z.
\end{equation}
\end{lem}

\proof One has $\Delta(1)=1\otimes 1$ and thus $
\Delta(nm)=n\otimes m \qqq n, m \in \Z$.
Taking $n=1$, \eqref{equhopf} gives $\varphi(m)=\varphi_2(m)$ for all $m\in \Z$ so that
$\varphi|_\Z=\varphi_2|_\Z$ and similarly taking $m=1$ we get
$\varphi|_\Z=\varphi_1|_\Z$.
\endproof

 \begin{lem} Let $\varphi\in\varphi_1\star_\Delta\varphi_2$ and let $J_j=Ker(\varphi_j)$ ($j=1,2$) be the associated ideals of $\cH=\Z[T]$. Then, one has
\begin{equation}\label{14h}
\varphi(x)=0\quad\forall x~\text{with}~ \Delta x\in J_1\otimes\cH+\cH\otimes J_2.
\end{equation}
\end{lem}
\proof Write a decomposition of the form
\[
\Delta x = \sum  x_i\otimes H_i + \sum  H_k'\otimes x_k', \quad x_i\in J_1,~ x_k'\in J_2.
\]
Then one has
\[
\varphi(x)\in \sum \varphi_1(x_i)\varphi_2(H_i)+\sum \varphi_1(H_k')\varphi_2(x_k')=0.
\]
\endproof
Let $J_j$ be ideals in $\cH$. The subset $J=J_1\otimes\cH+\cH\otimes J_2$ of $\cH\otimes\cH$ is an ideal of $\cH\otimes \cH$ and we set
\begin{equation}\label{15h}
J_1\star_\Delta J_2 =\{x\in\cH| \Delta x\in   J\}.
\end{equation}
Since $\Delta$ is a ring homomorphism, $J_1\star_\Delta J_2$ is an ideal of $\cH$.

\begin{lem}\label{sml} For any $\varphi\in \varphi_1\star_\Delta\varphi_2$, one has $Ker(\varphi_1)\star_\Delta Ker(\varphi_2)\subset Ker(\varphi)$.
\end{lem}
\proof Let $J_j=Ker(\varphi_j)$ ($j=1,2$).
For $\Delta x\in  J=J_1\star_\Delta J_2$  it follows  from \eqref{14h} that $\varphi(x)=0$.
\endproof

 \subsection{Hyperoperations on functions}\label{coprosect}

The above results
allow  one to define  the algebraic structure on functions, \ie on the elements of $\cD(R)=\Hom(\Z[T],R)$.
Here we use  the $2$ coproducts on $\cH=\Z[T]$ which are  uniquely defined  by
\begin{equation}\label{3h}
\Delta^+(T)=T\otimes 1+1\otimes T\in \cH\otimes \cH
\end{equation}
and
\begin{equation}\label{4h}
\Delta^\times(T)=T\otimes T\in \cH\otimes\cH.
\end{equation}
\begin{defn} \label{defcop} Let $R$ be a hyperring  and $\varphi_j\in Hom(\Z[T],R)$ be two functions. Let $\Delta^*$ be either $\Delta^+$ or $\Delta^\times$. One sets
\begin{equation}\label{defncop}
\varphi_1\star_{\Delta^*}\varphi_2 = \{\varphi\in Hom (\Z[T],R)|\varphi(x)\in \sum\varphi_1(x_{(1)})\varphi_2(x_{(2)})\}
\end{equation}
 for any decomposition $\Delta^*x = \sum x_{(1)}\otimes x_{(2)}$.
\end{defn}
We now give a general construction of the elements $\varphi\in\varphi_1\star_\Delta\varphi_2$.
\begin{lem}\label{tmlprel} Let $A$ be a ring and let $G\subset A^\times$  be a subgroup of the units of $A$. We denote by $R=A/G$ the quotient hyperring and we let $\epsilon: A \to A/G$  the projection map. For  $a\in A$, let $\tilde\varphi_a: \Z[T]\to A$ be the ring homomorphism given by $\tilde\varphi_a(p(T))= p(a)$ and we let
\begin{equation}\label{6h}
\varphi_a: \Z[T]\to A/G\qquad \varphi_a=\epsilon\circ\tilde\varphi_a
\end{equation}
to be the composite homomorphism (of hyperrings). Let $\Delta^*: \Z[T]\to \Z[T]\otimes_\Z \Z[T]$ be a coproduct and  $m: A\otimes A \to A$  the multiplication map. For $a,b\in A$, we let $c\in A$ be
\begin{equation}\label{7h}
c=\psi(T)\,, \ \psi=m\circ(\tilde\varphi_a\otimes\tilde\varphi_b)\circ\Delta^*: \Z[T]\to A.
\end{equation}
Then, for any decomposition $\Delta^* p = \sum p_{(1)}\otimes p_{(2)}$, one has
\begin{equation}\label{8h}
\varphi_c(p(T))\in \sum\varphi_a(p_{(1)})\varphi_b(p_{(2)}).
\end{equation}
\end{lem}
\proof  Notice that  $\psi$ is determined by  $\psi(T)=c$, thus  $\psi=\tilde\varphi_c$ and
 $\varphi_c=\epsilon\circ\psi$.
If $\Delta^* p = \sum p_{(1)}\otimes p_{(2)}$ one gets, by \eqref{7h}
\begin{equation}\label{10h}
\psi(p)=\sum \tilde\varphi_a( p_{(1)}) \tilde\varphi_b(p_{(2)})\in A.
\end{equation}
Since $\epsilon(\sum x_i)\in \sum\epsilon(x_i)$ for any $x_i\in A$, one obtains
\begin{equation}\label{11h}
\varphi_c(p)\in\sum \varphi_a(p_{(1)})\varphi_b(p_{(2)})
\end{equation}
as required.
\endproof

\begin{lem}\label{tml} With the notations of Lemma \ref{tmlprel} and  for any $a,b\in A$, one has
\begin{equation}\label{addi}
  \varphi_{a+b}\in  \varphi_a\star_{\Delta^+}\varphi_b
\end{equation}
\begin{equation}\label{multi}
 \varphi_{ab}\in  \varphi_a\star_{\Delta^\times}\varphi_b.
\end{equation}
\end{lem}
\proof
By applying Lemma \ref{tmlprel} to the coproduct $\Delta^*=\Delta^+$ {\it resp.} $\Delta^*=\Delta^\times$  , one gets $c=a+b$, resp. $c=ab$.\endproof

\section{Functions on $\Spec(\kras)$}\label{sectspeckras}

In this section we study the set of functions on $\Spec(\kras)$ with their hyperoperations.
The  map $\pi:\Hom(\Z[T],\kras)\to \Hom(\Z,\kras)$ (\cf \eqref{pi}) maps these functions to $\Spec\Z$ and we know (\cf Lemma \ref{restoz}) that the hyperoperations in $\cD(\kras)=\Hom(\Z[T],\kras)$  occur fiberwise. We shall thus describe separately these hyperoperations within elements of the same fiber $\pi^{-1}(p)$,  $p\in \Spec \Z$ (here, $\pi^{-1}(p)$ refers to $\pi^{-1}(\alpha)$, for $\alpha:\Z \to \kras$, \ie  we identify $p\in\Sp(\Z)$ with the kernel of a homomorphism $\alpha$).

\subsection{Functions on $\Sp(\kras)$: the fiber over $\{0\}$}\label{functions0} In the following, we consider the fiber of $\pi$ over the generic point $\{0\}$ of $\Spec\Z$. Let $\varphi\in \Hom(\Z[T],\kras)$, then $\varphi\in \pi^{-1}(\{0\})$ if and only if $\varphi(n)=1$ for all $n\in \Z$, $n\neq 0$. In turn, this holds if and only if $\varphi$ is the restriction to $\Z[T]$ of an element of $\Hom(\Q[T],\kras)$.
 This means that when $\varphi\in \pi^{-1}(\{0\})$,  one can extend  $\varphi$ to a homomorphism $\tilde\varphi: \Q[T]\to \kras$ by setting
\begin{equation}\label{16h}
\tilde\varphi(p(T))=\varphi(np(T))\qquad \forall n\neq 0,~~ np\in \Z[T].
\end{equation}
By taking a common multiple of the denominators appearing in the coefficients of $p(T)\in\Q[T]$,  one sees that  the definition of $\tilde\varphi$ is independent  of $n$ and determines a multiplicative map.  One also has
\[
\tilde\varphi(p_1+p_2)=\varphi(np_1+np_2)\in \varphi(np_1)+\varphi(np_2)=\tilde\varphi(p_1)+\tilde\varphi(p_2).
\]
Thus, one obtains the identification
\begin{equation}\label{zq}
   \pi^{-1}(\{0\})=\Hom(\Q[T],\kras)=\Spec\Q[T].
\end{equation}
Moreover, if $\Delta p = \sum p_{(1)}\otimes p_{(2)}$ holds  in the extension  of the coproduct to $\Q[T]\to \Q[T]\otimes \Q[T]$, one can find  $n_j\in\Z\setminus\{0\}$ ($j=1,2$) such that the equality
\[
\Delta(n_1n_2p) = \sum n_1p_{(1)}\otimes n_2 p_{(2)}
\]
only involves  elements of $\Z[T]$. This shows that one can set $\varphi_1\star_{\Delta^*}\varphi_2$ as in   \eqref{defncop} by implementing $\Q[T]$ rather than $\Z[T]$ in Definition \ref{defcop}.

\subsubsection{Hyperoperations on non-generic points}

We denote by $\delta$ the generic point of $\Spec\Q[T]=\pi^{-1}(\{0\})$ and we first determine the two hyperoperations on the complement $X=\pi^{-1}(\{0\})\setminus \{\delta\}$ of $\delta$ in $\Spec\Q[T]$. This complement  is the set of non-zero prime ideals of $\Q[T]$. We identify $X$ with the quotient of the field of the algebraic numbers $\bar\Q\subset \C$ by the action of the Galois group $\Aut_\Q(\bar \Q)$:
\begin{equation}\label{comp}
  X= \pi^{-1}(\{0\})\setminus \{\delta\}\simeq \bar \Q/\Aut(\bar \Q).
\end{equation}
To a non-zero prime ideal $\ffp$ of $\Q[T]$ one associates the roots in $\bar \Q$ of a generator of $\ffp$. These roots form an orbit for  the action of $\Aut(\bar \Q)$.
\begin{thm}\label{thmalg} The hyperoperations   $\varphi_1\star_{\Delta^*}\varphi_2$ ($*=+,\times$) of sum and product on $X= \pi^{-1}(\{0\})\setminus \{\delta\}$
coincide with the hyper-addition and hyper-multiplication on the hyperstructure $\bar \Q/\Aut(\bar \Q)$.
\end{thm}
\proof Let $\varphi_j\in\Hom(\Z[T],\kras)$ ($j=1,2$) be associated  to the prime  polynomials  $q_j\in \Z[T]$, \ie $Ker(\tilde\varphi_j)\subset \Q[T]$ is generated  by  $q_j$.  Let $\alpha\in Z(q_1)$ (resp. $\beta\in Z(q_2)$) be a zero of $q_1$ (resp. $q_2$). By Lemma~\ref{tml} (here applied to $A=\bar\Q$ and  $G = \bar\Q^\times$), the homomorphism $\varphi\in\Hom(\Z[T],\kras)$ defined by
$$
\varphi(P(T))=P(\alpha+\beta)G\in \bar\Q/\bar\Q^\times\simeq\kras
$$
belongs to $\varphi_1\star_{\Delta^+}\varphi_2$. A similar result holds for $\alpha\beta$ and $\Delta^\times$. It remains to show that all elements of $\varphi_1\star_{\Delta^*}\varphi_2$ are of this type. We first consider the addition, \ie  the set $\varphi_1\star_{\Delta^+}\varphi_2$.
Let $q\in\Q[T]$ be the monic polynomial  whose roots (in $\bar\Q$) are all the $\alpha+\beta$, with $\alpha\in Z(q_1)$ and $\beta\in Z(q_2)$, \ie $q(T)=\prod_{\alpha\in Z(q_1),\ \beta\in Z(q_2)} (T-\alpha-\beta)$. Let $n$ a non-zero integer such that $nq(T)\in \Z[T]$. Then, it is enough to show that $nq(T)\in Ker(\varphi)$ for any $\varphi\in \varphi_1\star_{\Delta^+}\varphi_2$. In fact, by applying Lemma \ref{sml} it is enough to show that
\begin{equation}\label{bel}
    nq(T)\in Ker(\varphi_1)\star_{\Delta^*} Ker(\varphi_2)
\end{equation}
for $\Delta^* = \Delta^+$. The polynomial $nq(T)$ is given, up to a multiplicative constant, by the resultant (\cf \cite{Bourbaki}, A IV, 6) of $q_1$ and $q_2$:
\begin{equation}\label{resu}
    q(Z)={\rm Resultant}_T(q_1(T),q_2(Z-T)).
\end{equation}
In particular, the polynomial $q(X+Y)$ is the resultant
\begin{equation}\label{resu1}
    q(X+Y)={\rm Resultant}_T(q_1(X-T),q_2(Y+T)).
\end{equation}
It follows from \cite{Bourbaki}, A IV, 6, Remark 4, by evaluation at $T=0$, that there exist polynomials
$A(X,Y)$ and $B(X,Y)$ with rational coefficients such that
$$
q(X+Y)=q_1(X)A(X,Y)+q_2(Y)B(X,Y).
$$
Then \eqref{bel} for $\Delta^*=\Delta^+$ easily follows.

One proceeds similarly with  the second co-product $\Delta^\times$. The polynomial $q$ is the  monic polynomial  whose roots (in $\bar\Q$) are the $\alpha\beta$ for $\alpha\in Z(q_1)$, $\beta\in Z(q_2)$ as above, \ie $q(T)=\prod_{\alpha\in Z(q_1),\ \beta\in Z(q_2)} (T-\alpha\beta)$. If $\beta=0$, then \eqref{4h} shows that the polynomial $q(T)=T$ belongs to $Ker(\varphi_1)\star_{\Delta^\times} Ker(\varphi_2)$. Thus we can assume that the roots of $q_2(T)$ are all non-zero. The polynomial $q$ is then, up to a multiplicative constant, the resultant
\begin{equation}\label{resu2}
    q(Z)={\rm Resultant}_T(q_1(T),T^mq_2(Z/T))
\end{equation}
where $m$ is the degree of $q_2$.
In particular, the polynomial $q(XY)$ is the resultant
\begin{equation}\label{resu3}
    q(XY)={\rm Resultant}_T(q_1(XT),T^mq_2(Y/T)).
\end{equation}
Then, it follows from \cite{Bourbaki}, A IV, 6, Remark 4, by evaluation at $T=1$, that there exist polynomials
$A(X,Y)$ and $B(X,Y)$ with rational coefficients such that
$$
q(XY)=q_1(X)A(X,Y)+q_2(Y)B(X,Y).
$$
This implies \eqref{bel} for $\Delta^*=\Delta^\times$.
\endproof

\subsubsection{Hyperaddition  with the generic point}

In this subsection we determine the hyperaddition law with the generic point $\delta\in \Spec(\Q[T])$. We shall make use, as before, of the identification  $X=\pi^{-1}(\{0\})\setminus \{\delta\}\simeq \bar \Q/\Aut(\bar \Q)$.
\begin{thm}\label{genericplus} In the fiber $\pi^{-1}(\{0\})=\Hom(\Q[T],\kras)=\Spec\Q[T]$, the hyperaddition with the generic point $\delta$ is given as follows\vspace{.05in}

$\bullet$~$\delta \star_{\Delta^+} \delta=\pi^{-1}(\{0\})$.\vspace{.05in}

$\bullet$~$\delta \star_{\Delta^+} \alpha=\pi^{-1}(\{0\})$,\quad $\forall\alpha\in\bar \Q\setminus \Q$.\vspace{.05in}

$\bullet$~$\delta \star_{\Delta^+} \alpha=\delta$,\quad$\forall\alpha\in  \Q $.
\end{thm}
\proof Let $A=\C$ and $G=\C^\times$. Then, with the notations of Lemma \ref{tmlprel}, one has
\begin{equation}\label{delta}
    \delta=\varphi_a\qqq a\in \C\setminus\bar\Q.
\end{equation}
Thus  Lemma \ref{tml} shows that $\delta \star_{\Delta^+} \delta=\pi^{-1}(\{0\})$ since any complex number can be written as the sum of two transcendental numbers. The same lemma also shows that $\delta \star_{\Delta^+} \alpha$ contains $\delta$ for any $\alpha\in\bar\Q$. The second and third equalities follow from Lemmas \ref{algplusgen1} and \ref{ratplusgen1} proven below.\endproof

\begin{lem}\label{algplusgen} Let $a(T)$ be an irreducible polynomial in $\Q[T]$
of degree $>1$ and let $J\subset \Q[T]$ be the prime ideal generated by  $a(T)$ in $\Q[T]$.
Let $P(T)$ be a  polynomial in $\Q[T]$. Let
\begin{equation}\label{dec}
    P(X+Y)=\sum A_j(X)B_j(Y)\,, \ A_j,B_j\in \Q[T]
\end{equation}
be a decomposition of $\Delta^+(P)$. Then, if the degree of $P$ is strictly positive
there exist at least two indices $j$ in the decomposition \eqref{dec} such that $A_j\neq 0$ and $B_j\notin J$.

If $P$ is a constant polynomial and $i=1$  is the only index in \eqref{dec} for which $A_i\neq 0$ and $B_i\notin J$, then $A_1(X)=d$ is a non-zero constant polynomial and $a(T)$ divides $d B_1(T)-c$.
\end{lem}

\proof Let assume first that $P$ is non-constant.
Also, let assume that only $B_1$ does not belong to the prime ideal $J$ generated by  $a(T)$ in $\Q[T]$. Let $\alpha_\ell\in \bar \Q$, ($\ell=1,2$) be two distinct roots of $a(T)$. Since $B_j(\alpha_\ell)=0$ for $j\neq 1$, one has, using \eqref{dec}
$$
P(X+\alpha_\ell)=A_1(X)B_1(\alpha_\ell)\,, \ \ell=1,2.
$$
 It follows that $B_1(\alpha_\ell)\neq 0$ and that
$$
P(X+\alpha_1)/B_1(\alpha_1)=P(X+\alpha_2)/B_1(\alpha_2)
$$
so that
$$
P(X+\alpha_1-\alpha_2)=\lambda P(X)\,, \ \ \lambda=B_1(\alpha_1)/B_1(\alpha_2).
$$
But since $\alpha_1-\alpha_2\neq 0$ and $\lambda\neq 0$ this yields infinitely many zeros for $P$, thus we derive a contradiction.

Assume now that  $P=c$ is constant and $i=1$  is the only index in \eqref{dec}  for which $A_i\neq 0$ and $B_i\notin J$. Let $\alpha\in \bar \Q$ be a root of $a(T)$. Then \eqref{dec} implies that
$$
c=A_1(X)B_1(\alpha).
$$
Thus $A_1(X)=d$ is a non-zero constant. Hence $B_1(\alpha)=c/d$ is independent of the choice of $\alpha$ and $d B_1(T)-c$ is divisible by $a(T)$.
 \endproof\vspace{.05in}

Let $a(T)$ be an irreducible polynomial in $\Q[T]$, we set
\begin{equation}\label{phiapsia}
    \psi_a=\varphi_\alpha\qqq \alpha\in \bar\Q\,, \ a(\alpha)=0.
\end{equation}

\begin{lem}\label{algplusgen1} Let $a(T)$ be an irreducible polynomial in $\Q[T]$
of degree $>1$. Then
\begin{equation}\label{genplusal}
  \delta\star_{\Delta^+}\psi_a=\pi^{-1}(0).
\end{equation}
\end{lem}

\proof Let $\varphi\in \pi^{-1}(0)$. In order to show that $\varphi\in\delta\star_{\Delta^+}\psi_a$, we need to prove that for any   polynomial  $P\in\Q[T]$
 and any decomposition as \eqref{dec} one has
\begin{equation}\label{pluscond}
\varphi(P)\in \sum\delta(A_j)\psi_a(B_j).
\end{equation}
If there are two indices $i$ for which $A_i\neq 0$, $B_i\notin J$, then \eqref{pluscond} follows from $\psi_a(B_i)=1$ and $\delta(A_i)=1$ since we then derive
\begin{equation}\label{overallsum}
\sum\delta(A_j)\psi_a(B_j)=\{0,1\}.
\end{equation}
Thus Lemma \ref{algplusgen} shows that \eqref{pluscond} holds when $P$ is non-constant.
If $P=c\neq 0$ is constant and  $i=1$  is the only index $i$ for which $A_i\neq 0$, $B_i\notin J$ then by applying Lemma \ref{algplusgen} we get that  $A_1(X)=d$ is a constant and $a(T)$ divides $d B_1(T)-c$, thus $\psi_a(B_1)= \psi_a(dB_1)=\psi_a(c)$ and both sides of \eqref{pluscond} are equal to $\{1\}$. Finally if $P=0$ then either $B_j\in J$ for all $j$ with $A_j\neq 0$ or there are two indices $j$ for which $A_j\neq 0$ and $B_j\notin J$. In both cases on has \eqref{pluscond}.
The  inclusion $ \delta\star_{\Delta^+}\psi_a\subset\pi^{-1}(0)$ follows from Lemma \ref{restoz}. \endproof

\begin{lem}\label{ratplusgen1} Let $a(T)$ be a polynomial in $\Q[T]$
of degree $1$. Then
$$
  \delta\star_{\Delta^+}\psi_a= \delta.
$$
\end{lem}

\proof Let $a(T)=mT-n$ where $m\neq 0$ and $m,n\in \Z$. Let $\varphi\in \delta\star_{\Delta^+}\psi_a$, we want to show that $\varphi(P)=1$ for any polynomial $P(T)\neq 0$. Using the Taylor expansion of $P$ at $Y=\frac n m$, we have
$$
P(X+Y)=P(X+\frac nm)+D(X,Y)(Y-\frac nm)
$$
with $D(X,Y)\in \Q[X,Y]$. By multiplying both sides of the above equality by a non-zero integer $k$ to get rid of the denominators, we obtain an equality of the form
$$
d P(X+Y)=A(X)+E(X,Y) a(Y)
$$
where $A$ and $E$ have both integral coefficients. Then, the definition of $\delta\star_{\Delta^+}\psi_a$ shows that
$$
\varphi(P)=\varphi(kP)\in \delta(A)=\{1\}.
$$
Thus $\varphi(P)=1$ for any polynomial $P(T)\neq 0$ and $\varphi=\delta$.
\endproof

\begin{rem}\label{remarkassoc}{\rm By Theorem \ref{thmalg}, the hyperaddition in $X=\pi^{-1}(\{0\})\setminus \{\delta\}$ defines a {\em canonical} hypergroup. Using Theorem \ref{genericplus}, one checks that the presence of the generic element $\delta$ does not spoil the associativity. Indeed, the sum $(x\star_{\Delta^+}y)\star_{\Delta^+}z$ of three elements one of which is $\delta$ is equal to $\pi^{-1}(0)$, unless the two remaining elements are in $\Q$, in which case the sum is equal to $\delta$. Note also that one has $x\star_{\Delta^+}y\subset \Q$ only if $x$ and $y$ are in $\Q$. However the reversibility property for hypergroups  no longer holds since for $\alpha\in \Q$ and $\beta\in \bar\Q\setminus \Q$ one has
$$
\alpha\in \delta-\beta\,, \ \beta\notin \delta -\alpha.
$$
}\end{rem}

\subsubsection{Hyper-multiplication with the generic point}

We shall keep using the identification \eqref{comp}. Note  that $0$ is an absorbing element for the hyperoperation  $\star_{\Delta^\times}$.
\begin{thm}\label{genericprod} In the fiber $\pi^{-1}(\{0\})=\Hom(\Q[T],\kras)=\Spec\Q[T]$, the hyper-multiplication with the generic point $\delta$ is given by\vspace{.05in}

$\bullet$~$\delta\star_{\Delta^\times}\delta=\pi^{-1}(0)\setminus\{0\}$.\vspace{.05in}

$\bullet$~$\delta \star_{\Delta^\times} \alpha=\pi^{-1}(0)\setminus\{0\}$,~ $\forall \alpha\in\bar \Q$~ $\alpha^n\notin \Q$, $\forall n>0$.\vspace{.05in}

$\bullet$~$\delta \star_{\Delta^\times} \alpha=\{\delta\}$,~ $\forall\alpha\in\bar \Q$, $\alpha\neq 0$,~ $\alpha^n\notin \Q$ for some $n>0$.
\end{thm}
\proof Let $A=\C$ and $G=\C^\times$. Then, with the notations of Lemma \ref{tmlprel}, one has
\eqref{delta}. Lemma \ref{tml} implies the first equality since any non-zero complex number can be written as the product of two transcendental numbers. The same lemma also shows that $\delta \star_{\Delta^\times} \alpha$ contains $\delta$ for any non-zero $\alpha$. The second and third equalities follow from Lemmas \ref{algtimesgen1} and \ref{algtimesgen2} here below.\endproof

\begin{lem}\label{algtimesgen} Let $a(T)$ be an irreducible polynomial in $\Q[T]$
which admits two non-zero roots $\alpha_1,\alpha_2\in\bar\Q$ whose ratio is not a root of unity and let $J\subset \Q[T]$ be the prime ideal generated by  $a(T)$ in $\Q[T]$.
Let $P(T)$ be a  polynomial in $\Q[T]$. Let
\begin{equation}\label{decprod}
    P(XY)=\sum A_j(X)B_j(Y)\,, \ A_j,B_j\in \Q[T]
\end{equation}
be a decomposition of $\Delta^\times(P)$. If  $P\neq 0$ has a non-zero root, then
there exist at least two indices $j$ in the decomposition \eqref{decprod} such that $A_j\neq 0$ and $B_j\notin J$.

If $P=cT^n$, $c\in\Q$  and $i=1$  is the only index  for which $A_i\neq 0$, $B_i\notin J$ in \eqref{decprod}, then $A_1(X)=dX^n$, with $d\in\Q$  and $a(T)$ divides $d B_1(T)-cT^n$.
\end{lem}

\proof Assume first that $P$ has a non-zero root. We also
assume that $A_1\neq 0$ and only $B_1$ does not belong to the prime ideal $J$ generated by  $a(T)$ in $\Q[T]$. Let $\alpha_\ell$, $\ell=1,2$, be two distinct non-zero roots of $a(T)$ whose ratio is not a root of unity. Since $B_j(\alpha_\ell)=0$ for $j\neq 1$, one has
$$
P(X\alpha_\ell)=A_1(X)B_1(\alpha_\ell)\,, \ \ell=1,2.
$$
 Then, it follows that $B_1(\alpha_\ell)\neq 0$ and
$$
P(X\alpha_1)/B_1(\alpha_1)=P(X\alpha_2)/B_1(\alpha_2)
$$
so that
$$
P(X\alpha_1/\alpha_2)=\lambda P(X)\,, \ \ \lambda=B_1(\alpha_1)/B_1(\alpha_2).
$$
But since $\alpha_1/\alpha_2$ is not a root of unity,  $\lambda\neq 0$ and $P$ has a non-zero root, this argument yields infinitely many zeros for $P$, thus it produces a contradiction.

Assume now that  $P=cT^n$, for $c\in\Q$  and that $i=1$  is the only index $i$ in the decomposition \eqref{decprod} of $P(XY)$ for which $A_i\neq 0$ and $B_i\notin J$. Let $\alpha\in \bar \Q$, $\alpha\neq 0$ be a root of $a(T)$. Then \eqref{decprod} implies that
$$
cX^n\alpha^n=A_1(X)B_1(\alpha).
$$
Thus $A_1(X)=dX^n$ where $d\neq 0$, $d\in\Q$. Hence $B_1(\alpha)=c\alpha^n/d$ and $d B_1(T)-cT^n$ is divisible by $a(T)$ since it vanishes on all roots of $a(T)$.  \endproof

\begin{lem}\label{algtimesgen1} Let $a(T)$ be an irreducible polynomial in $\Q[T]$
which admits two non-zero roots $\alpha_1,\alpha_2\in\bar\Q$ whose ratio is not a root of unity. Then (\cf\eqref{phiapsia})
\begin{equation}\label{gentimesal}
  \delta\star_{\Delta^\times}\psi_a=\pi^{-1}(0)\setminus\{0\}.
\end{equation}
\end{lem}
\proof Let $\varphi\in \pi^{-1}(0)$ with $\varphi\neq \psi_T$.
We show that for any $P\in \Q[T]$ and any decomposition \eqref{decprod}, one has
 \begin{equation}\label{timescond}
\varphi(P)\in \sum\delta(A_j)\psi_a(B_j).
\end{equation}
We first assume that $P$  admits a non-zero root. By Lemma \ref{algtimesgen} there are at least two indices $i=1,2$ for which $A_i\neq 0$, $B_i\notin J$. Thus $\delta(A_i)=1$,
$\psi_a(B_i)=1$ and so
\begin{equation}\label{overallsum1}
\sum\delta(A_j)\psi_a(B_j)=\{0,1\}.
\end{equation}
The remaining polynomials  $P(T)$  have only the zero root $T=0$ and  thus they are of the form $cT^n$ where $c\in \Q$. For $c=0$, Lemma \ref{algtimesgen} shows that either all $B_i\in J$ or there are at least two indices $i=1,2$ for which $A_i\neq 0$, $B_i\notin J$. Thus in this case \eqref{timescond} holds. We can thus assume $c\neq 0$.     Since $\varphi\neq \psi_T$ one has $\varphi(P)=1$. We prove that the right hand side of \eqref{timescond} always contains $1$. This statement holds when there are at least two indices $i$ for which $A_i\neq 0$, $B_i\notin J$.
By Lemma \ref{algtimesgen}, if  $i=1$  is the only index $i$ for which $A_i\neq 0$, $B_i\notin J$ then $A_1(X)=dX^n$, where $d\neq 0$ and $a(T)$ divides $d B_1(T)-cT^n$. It is enough to show that $\psi_a(B_1)=1$. From the inclusion $\psi_a(x+y)\subset \psi_a(x)+\psi_a(y)$ and the equality
 $$
 \psi_a(d B_1(T)-cT^n)=0
 $$
it follows
$$
\psi_a(B_1)=\psi_a(dB_1)=\psi_a(cT^n)=1.
$$
 Finally, note that $0$ does not belong to $\delta\star_{\Delta^\times}\psi_a$ since $\psi_T(T)=0$ while $\Delta^\times(T)=T\otimes T$ and $\delta(T)=\psi_a(T)=1$.
 \endproof

 \begin{lem}\label{algtimesgen1} Let $a(T)$ be an irreducible polynomial in $\Q[T]$ not proportional to $T$,
such that the ratio of any two  roots $\alpha_1,\alpha_2\in\bar\Q$ of $a(T)$  is  a root of unity.
Let $P(T)$ be an irreducible  polynomial in $\Q[T]$. Then $P(T)$ has a multiple $M(T)=C(T)P(T)$ such that $\Delta^\times(M)$ admits a decomposition
\begin{equation}\label{decprodbis}
    M(XY)= A(X)B(Y)+D(X,Y)a(Y)
\end{equation}
where $A\neq 0$ and $B$ is not divisible by $a$.
\end{lem}

\proof One has $(\alpha_1/\alpha_2)^n=1$ for a suitable $n$ and all roots of $a(T)$. Thus $\alpha^n=\lambda\in \Q$ for all roots $\alpha$ of $a$. This shows that $a(T)$ is an irreducible factor of $T^n-\lambda$. Let $N\in \Q[T]$ and consider the polynomial $H(T)=N(T^n)$.
Next we show that one has a decomposition
\begin{equation}\label{decN}
  H(XY)=A(X)+D(X,Y)a(Y).
\end{equation}
We let $A(X)=N(X^n\lambda)$ and note that $K(X,Y)=H(XY)-N(X^n\lambda)$ vanishes when $Y^n=\lambda$ and hence when $Y$ is a root of $a(Y)$. Thus $K(X,Y)$ is divisible by $a(Y)$
and it can be written in the form $D(X,Y)a(Y)$; this proves \eqref{decN}.
It remains to show that any irreducible polynomial $P(T)$ has a multiple $H(T)=C(T)P(T)$ of the form $H(T)=N(T^n)$. Let
$$
P(T)=c\prod_k (T-\alpha_k)
$$
be the factorization of $P(T)$ in $\bar\Q[T]$.
Let
$$
H(T)=N(T^n)\,, \ \ N(T)=c\prod_k (T-\alpha_k^n).
$$
Then $N(T)\in \Q[T]$ and $P$ divides $H$. This provides the required conclusion.
\endproof

\begin{lem}\label{algtimesgen2} Let $a(T)$ be an irreducible polynomial in $\Q[T]$ not proportional to $T$,
such that the ratio of any two  roots $\alpha_1,\alpha_2$ of $a(T)$  is  a root of unity. Then
\begin{equation}\label{gentimesal}
  \delta\star_{\Delta^\times}\psi_a= \delta.
\end{equation}
\end{lem}

\proof Let $\varphi \in   \delta\star_{\Delta^\times}\psi_a$. Let $P(T)$ be an irreducible  polynomial in $\Q[T]$. We want to show that $\varphi(P)=1$. By Lemma \ref{algtimesgen1},  $P(T)$ has a multiple $M(T)=C(T)P(T)$ fulfilling \eqref{decprodbis}.
Then
$$
\varphi(M)\in \delta(A)\psi_a(B)=1
$$
thus $\varphi(M)=1$ and $\varphi(P)=1$. This shows that only $\delta$ may belong to $\delta\star_{\Delta^\times}\psi_a$ and this is the case as  follows from the proof of Theorem \ref{genericprod}.
\endproof

\begin{rem}\label{remarkassocbis}{\rm By Theorem \ref{thmalg}, the hypermultiplication in $\pi^{-1}(\{0\})\setminus \{0,\delta\}$ defines a {\em canonical} hypergroup. Using Theorem \ref{genericprod} one checks that the presence of the generic element $\delta$ does not spoil the associativity. Indeed the product $(x\star_{\Delta^\times}y)\star_{\Delta^\times}z$ of three non-zero elements one of which is $\delta$ is equal to $\pi^{-1}(0)\setminus \{0\}$ unless the two remaining elements are in $\Q^{\rm root}=\{\alpha\mid\exists n,\alpha^n\in\Q\}$ and in that case the product is equal to $\delta$. Note that one has $x\star_{\Delta^\times}y\subset \Q^{\rm root}$ only if $x$ and $y$ are in $\Q^{\rm root}$. However the reversibility property for hypergroups  no longer holds.
}\end{rem}

\subsection{Functions on $\Spec(\kras)$: the fiber over $\{p\}$}\label{functionsp}

Let $p$ be a prime integer. In this section we compute the hyper-addition and the hyper-multiplication in the fiber $\pi^{-1}(p)$ of $\pi: \Hom(\Z[T],\kras)\to\Hom(\Z,\kras)$. We let $\Omega$ be an algebraic closure of the field of fractions $\F_p(T)$.
\begin{thm}\label{thmbijp} The following map determines a bijection of $\Omega/\Aut(\Omega)$
with the fiber $\pi^{-1}(p)\subset \Hom(\Z[T],\kras)$
\begin{equation}\label{bijp}
    \alpha\in \Omega\mapsto \varphi_\alpha\,, \ \ \varphi_\alpha(P(T))=P(\alpha)\Omega^\times\in \Omega/\Omega^\times\simeq\kras.
\end{equation}
The hyperoperations $\varphi_1\star_{\Delta^*}\varphi_2$ ($*=+,\times$) on  the fiber $\pi^{-1}(p)$
coincide with the hyper-addition and hyper-multiplication of the hyperstructure $\Omega/\Aut(\Omega)$.
\end{thm}
\proof For $\alpha\in \Omega$ one has by construction $\varphi_\alpha\in \pi^{-1}(p)\subset \Hom(\Z[T],\kras)$. Conversely, let $\varphi\in \pi^{-1}(p)$, then the kernel of $\varphi$ determines a prime ideal $J\subset\F_p[T]$. If $J=\{0\}$, then $\varphi=\varphi_\alpha$ for any $\alpha\in \Omega\setminus \bar\F_p$. By \cite{Bourbaki} Proposition 9 (Chapter V, \S XIV, AV 112),
the group $\Aut(\Omega)$ acts transitively on the complement  $\Omega\setminus\bar\F_p$
of the algebraic closure of $\F_p$. If  $J\neq\{0\}$, then it  is generated by an irreducible (separable) polynomial $a(T)\in \F_p[T]$ and $\varphi=\varphi_\alpha$ for any root of $a(T)$ in $\bar\F_p\subset \Omega$. The set of roots of $a(T)$ in $\bar\F_p\subset \Omega$ is a single orbit of the action of $\Aut(\bar\F_p)$ and hence of $\Aut(\Omega)$
(\cf \cite{Bourbaki} Corollary 1, Chapter V, \S XIV, AV 111). This proves the first statement.
The proof of Theorem \ref{thmalg}, with $\F_p[T]$ in place of $\Q[T]$ shows that the hyper-operations   $\varphi_1\star_{\Delta^*}\varphi_2$ ($*=+,\times$) of sum and product on the complement of the generic point in $\pi^{-1}(p)$
coincide with the hyper-addition and hyper-multiplication on the hyperstructure $\bar \F_p/\Aut(\bar \F_p)$. It remains to determine these operations when the generic point $\delta_p$ is involved. It  follows from Lemma \ref{tml} applied to  $\kras\simeq\Omega/\Omega^\times$ that
$$
\delta_p\star_{\Delta^+}\delta_p= \pi^{-1}(p).
$$
The end of the proof then follows from  Lemma \ref{algplusgenp1} and  Lemma \ref{algmultgenp1}. \endproof

We now use the existence of enough ``additive" polynomials in characteristic $p$ (\cf \cite{goss}).

\begin{lem}\label{algplusgenp} Let $a(T)\in \F_p[T]$ be an irreducible polynomial of degree $n>0$. Let $S_m(T)=T^{p^m}-T$, with $n|m$. Then, there exists $B(T)\in \F_p[T]$ such that
\begin{equation}\label{decp}
    S_m(Y)=B(Y)a(Y)
\end{equation}
\begin{equation}\label{decpbis}
    S_m(X+Y)=S_m(X)+B(Y)a(Y).
\end{equation}
\end{lem}

\proof Every root of $a(T)$ in $\bar\F_p$ is a root of $S_m(T)$, thus $a$ divides $S_m$.  The second statement follows from the equality $S_m(X+Y)=S_m(X)+S_m(Y)$. \endproof

\begin{lem}\label{algplusgenp1} Let $a(T)\in \F_p[T]$ be an irreducible polynomial of degree $n>0$. Then
$$
  \delta_p\star_{\Delta^+}\psi_a= \delta_p.
$$
\end{lem}

\proof Let $P(T)\in \F_p[T]$ be an irreducible polynomial of degree $k$ and let $\varphi \in \delta_p\star_{\Delta^+}\psi_a$. Then we show that $\varphi(P)=1$. By Lemma \ref{algplusgenp}, $S_m$ is a multiple of $P$ for any integer multiple $m$ of the degree $k$ of $P$. By taking $m=kn$ and using Lemma \ref{algplusgenp}, we obtain
$$
\Delta^+(S_m)=S_m\otimes 1+1\otimes Ba
,\ \
\varphi(S_m)\in \delta_p(S_m)\psi_a(1)+\delta_p(1)\psi_a(Ba)=\{1\}.
$$
Since $P$ divides $S_m$, one gets $\varphi(P)=1$. This shows that only $\delta_p$ can belong to $\delta_p\star_{\Delta^+}\psi_a$. Finally we prove that $\delta_p\star_{\Delta^+}\psi_a\ni \delta_p$. This follows from Lemma \ref{tml} applied to  $\Omega/\Omega^\times\simeq\kras$, since $\delta_p=\varphi_\gamma$ for any $\gamma\in \Omega\setminus \bar\F_p$ and $\psi_a=\varphi_\alpha$ for any root $\alpha$ of $a$. \endproof

\begin{lem}\label{algmultgenp} Let $a(T)\in \F_p[T]$ be an irreducible polynomial of degree $n>0$, not proportional to $T$. Let $U_k(T)=T^{p^k-1}-1$, with $n|k$. Then, there exists $B(T)\in \F_p[T]$ such that
\begin{equation}\label{decp1}
    U_k(Y)=B(Y)a(Y)
\end{equation}
\begin{equation}\label{decp1bis}
    U_k(XY)=U_k(X)Y^{p^k-1}+B(Y)a(Y).
\end{equation}
\end{lem}

\proof Every root of $a(T)$ in $\bar\F_p$ is a root of $U_k$ and thus $a$ divides $U_k$.  The second statement follows from the equality
$$
(XY)^{p^k-1}-1=(X^{p^k-1}-1)Y^{p^k-1}+(Y^{p^k-1}-1).
$$
 \endproof

\begin{lem}\label{algmultgenp1} Let $a(T)\in \F_p[T]$ be an irreducible polynomial of degree $n>0$, not proportional to $T$. Then
$$
  \delta_p\star_{\Delta^\times}\psi_a= \delta_p.
$$
\end{lem}
\proof
Let $P(T)\in \F_p[T]$ be an irreducible polynomial not proportional to $T$. Then, by taking for $k$ a common multiple of $n$ and of the degree of $P$, one gets from \eqref{decp} that $U_k(T)$ is a multiple of $P(T)$ and it fulfills \eqref{decp1bis}. Let $\varphi \in \delta_p\star_{\Delta^\times}\varphi_a$. We show that $\varphi(P)=1$. It is enough to prove that $\varphi(U_k)=1$. From \eqref{decp1bis}, it follows that
$$
\Delta^\times(U_k)=U_k\otimes T^{p^k-1}+1\otimes Ba, \ \ \varphi(U_k)\in \delta_p(U_k)\psi_a(T^{p^k-1})+\delta_p(1)\psi_a(Ba)=1.
$$
One has $\Delta^\times(T)=T\otimes T$ and thus $\varphi(T)\in \delta_p(T)\psi_a(T)=\{1\}$ so that $\varphi(T)=1$ and one gets $\varphi=\delta_p$.
Since Lemma \ref{tml} shows that $\delta_p\in \delta_p\star_{\Delta^\times}\psi_a$, one gets the required conclusion.
\endproof

\section{Functions on $\Sp(\sign)$}\label{sectspecsign}

We recall that $\sign = \{\pm 1, 0\}$ is the hyperfield of signs. We start by describing the functor $\Hom(\cdot,\sign)$ on ordinary rings. We let
\begin{equation}\label{1k}
|~|: \sign \to \kras
\end{equation}
be the homomorphism absolute value. Thus, for a given (commutative) ring $A$ and an element $\varphi\in\Hom(A,\sign)$,  one has: $|~|\circ\varphi = |\varphi| \in \Hom(A,\kras)$. This composite map is determined by its kernel which is a prime ideal of $A$. We recall (\cf\cite{wagner} Proposition 2.11) the following result
\begin{prop} \label{order} An element $\varphi\in\Hom(A,\sign)$ is determined by\vspace{.05in}

$a)$~ its kernel $Ker(\varphi)\in\Sp(A)$\vspace{.05in}

$b)$~ a total order on the field of fractions of $A/\wp$, $\wp= Ker(\varphi)$.
\end{prop}

Note that given  a prime ideal $\wp\subset A$ and a total order on the field $F$ of fractions of $A/\wp$, the corresponding homomorphism $\varphi\in\Hom(A,\sign)$ is the composite
\begin{equation}\label{2k}
\varphi: A \to A/\wp \to F \to F/F_+^\times \simeq \sign.
\end{equation}
The kernel of  $\varphi$ is $\wp$ and $\varphi^{-1}(1)$ determines the order on $F$.

\subsection{Description of $\Hom(\Z[T],\sign)$}

Notice that since $1+1=1$ in $\sign$, the set $\Hom(\Z,\sign)$ contains only one element which corresponds to $Ker(\varphi) = (0)$ (here $A=\Z$) and the usual order on $\Q$. We call this element $\varphi$ the {\em sign}. Thus, for $\varphi\in\Hom(\Z[T],\sign)=\cD(\sign)$ a function over $\Sp(\sign)$, the restriction of $\varphi$ to $\Z$ is equal to the sign and there is a unique extension $\tilde\varphi$ of $\varphi$ to $\Hom(\Q[T],\sign)$ given by
\begin{equation}\label{3k}
\tilde\varphi(p(T)) = \varphi(np(T)),\quad n>0,~np(T)\in\Z[T].
\end{equation}
This construction determines  an identification $\Hom(\Z[T],\sign)\simeq\Hom(\Q[T],\sign)$.  The functions on $\Sp(\sign)$  are thus determined by elements of  $\Hom(\Q[T],\sign)$, \ie by\vspace{.05in}

$a)$~ a prime ideal $J\subset\Q[T]$\vspace{.05in}

$b)$~a total order on the field of fractions of $\Q[T]/J$.\vspace{.1in}

The description of the set $\Hom(\Q[T],\sign)$ is given in \cite{wagner} Proposition 2.12.

\begin{prop}\label{mprop}
The elements of $\cD(\sign)=\Hom(\Z[T],\sign)$ are described by
\begin{equation}\label{omega}
    \omega_\lambda(P(T))={\rm Sign}(P(\lambda))\qqq \lambda \in [-\infty,\infty]
\end{equation}
and, for $\lambda\in \bar\Q\cap \R$, by the two elements
\begin{equation}\label{omegapm}
\omega_\lambda^\pm(P(T))=\lim_{\epsilon\to 0+}{\rm Sign}(P(\lambda\pm \epsilon)).
\end{equation}
\end{prop}

One has the natural map
\begin{equation}\label{remap}
Re: \Hom(\Q[T],\sign)\to [-\infty,\infty]
\end{equation}
defined by
\begin{equation}\label{4k}
Re(\varphi) = \sup\{a\in\Q, \varphi(T-a)=1\}\in[-\infty,\infty].
\end{equation}
The set of rational numbers involved in \eqref{4k} is a Dedekind cut (with a possible largest element). One thus gets
\begin{equation}\label{cut}
   \varphi(T-a) = 1, \ \forall a<Re(\varphi)\,, \ \varphi(T-a) = -1, \ \forall a>Re(\varphi).
\end{equation}

\begin{cor}\label{sign}

 $(1)$~Let $\varphi\in\Hom(\Q[T],\sign)$. If $Ker(\varphi)\neq (0)$ is generated by the prime polynomial $q(T)$, then $Re(\varphi)=\alpha\in\R$ is a root of $q(\alpha)=0$, and $\varphi=\omega_\alpha$.\vspace{.05in}

   $(2)$~There is only one $\varphi\in\Hom(\Q[T],\sign)$ such that $Re(\varphi)=+\infty$ (\resp $-\infty$) and it is given by
\begin{equation}\label{9k}
\varphi(p(T)) = \lim_{\alpha\to\infty}{\rm Sign}(p(\alpha))\quad (\resp \lim_{\alpha\to-\infty}{\rm Sign}(p(\alpha))).
\end{equation}

   $(3)$~If $\lambda\in\R\setminus\bar\Q$, then $\omega_\lambda$ is the only element $\psi$ such that $Re(\psi) = \lambda$.\vspace{.05in}

   $(4)$~If $Re(\varphi) = \alpha\in\bar\Q$ and $Ker(\varphi)= (0)$, then $\varphi=\omega_\alpha^s$ where $s=\varphi(q)\in \{\pm 1\}$ and $q\in \Q[T]$ is an irreducible polynomial such that $q(\alpha)=0$ and $q'(\alpha)>0$.
 \end{cor}
\proof All statements are straightforward consequences of the concrete description  of the elements of $\Hom(\Q[T],\sign)$ given in Proposition \ref{mprop}.
Next, we shall  explain briefly how the above statements can be proven directly using \cite{Bourbaki} Chapter VI.

$(1)$~When $Ker(\varphi)\neq (0)$, the quotient $F=\Q[T]/Ker(\varphi)$ is a finite extension of $\Q$ and a total order on $F$ is necessarily archimedean and hence it is given by an order embedding $F\subset \R$ \ie the map $P(T)\mapsto P(\alpha)$, where $\alpha$ is a real root of $q(T)$. One then gets
$\varphi=\omega_\alpha$ and $Re(\varphi)=\alpha$.

$(2)$~If $Re(\varphi)=+\infty$, one has $T\geq n$ for all $n\in \N$ for the corresponding order on $\Q[T]$ and the proof of Proposition 4 of \cite{Bourbaki} A.VI.24 shows that \eqref{9k} holds.

$(3),(4)$~Assume that $Re(\varphi) = \alpha\in\R$ and $Ker(\varphi)= (0)$. Then, by Proposition \ref{mprop} one gets a total order on the field $\Q(T)$ of rational fractions and  for some $n\in\N$ one has $-n<T<n$. Thus, all polynomials $p(T)$ are finite, \ie they belong to an interval $[-m,m]$ for some finite $m$. If all fractions $p(T)/q(T)$ are finite, then the order on $\Q(T)$ is archimedean and thus one derives an order embedding $\Q(T)\to \R$ which shows that $\varphi = \omega_\alpha$, $\alpha=Re(\varphi)$. Since $p(T)$ is finite, this can fail  only when some $q(T)\neq 0$ is infinitesimal
\begin{equation}\label{11k}
\pm q(T)<1/n\qquad \forall n\in\N.
\end{equation}
Thus the statements $(3),(4)$ follow from the next lemma.\endproof

\begin{lem} Let $\leq$ be a total order on $\Q[T]$ such that $T$ is finite. Let $J$ be the set of polynomials $q(T)\in\Q[T]$  fulfilling \eqref{11k}. Then\vspace{.05in}

$\bullet$~$J$ is a prime ideal\vspace{.05in}

$\bullet$~ The sign of $p\notin J$ only depends upon the class of $p$ in the quotient $\Q[T]/J$ and there exists a real root $\alpha$ of the monic polynomial $a(T)$ generating $J$ such that
      \begin{equation}\label{palpha}
        p>0\iff p(\alpha)>0\qqq p\notin J.
      \end{equation}
 $\bullet$~The total order $\leq$ is uniquely determined by $\alpha$ and the sign of $a(T)$.
\end{lem}
\proof Since $p(T)$ is finite, $\forall p(T)$, $J$ is an ideal. Moreover, if $q_j(T)>1/n_j$ one has $q_1(T)q_2(T)>1/n_1n_2$ and one derives that $J^c=\Q[T]\setminus J$ is multiplicative, thus $J$ is prime.

Next we show that the total order of $\Q[T]$ defines a total ordering on the quotient $\Q[T]/J$. This follows by noticing that if $p(T)>0$ in $\Q[T]$ and $p(T)\notin J$, then $p(T)+r(T)>0$ for any $r(T)\in J$. Indeed, since $p\notin J$, one has $p(T)>1/n$ for some $n\in\N$, but then since $-r<1/n$, one derives $p+r>0$. Thus, one gets a total order on $\Q[T]/J$. As in the proof of $(1)$ in Corollary \ref{sign}, this order is archimedean and thus it is produced by an algebraic number $\alpha$, root of the generator $a(T)$ of the prime ideal $J$. Thus \eqref{palpha} follows.
 Let $p(T)\in\Q[T]$, then there exists a unique $n$ such that $p\in J^n$, $p\notin J^{n+1}$. On has $p = a^n A$ where $A\notin J$. Since $A\notin J$, the sign of $A$ is given by the sign of $A(\alpha)\neq 0$. Thus the sign of $p$ is determined by the sign of $a$. \endproof

\subsection{Operations on functions on $\Spec(\sign)$}

We now investigate the hyper-structure on the set of functions on $\Sp(\sign)$. We consider the {\em finite} functions, \ie the elements of
\begin{equation}\label{finitefunct}
\cD_{\rm finite}(\sign)=\{\varphi\in\Hom(\Q[T],\sign)\mid Re(\varphi)\in \R\}.
\end{equation}
Let $F = \R(\epsilon)$ be the ordered field of rational fractions in $\epsilon$, with $\epsilon>0$ and infinitesimal. In particular, the sign of a polynomial $p(\epsilon) = a_k\epsilon^k+a_{k+1}\epsilon^{k+1}+\cdots+a_n\epsilon^n$ is given by the $\ss(a_k)$ \ie by the limit
\begin{equation}\label{8k}
\lim_{\epsilon\to 0+}\ss(p(\epsilon)).
\end{equation}
Any finite function  $\varphi$ with $Re(\varphi)=\alpha$ can be obtained as the composite
\begin{equation}\label{19k}
\Q[T]\stackrel{\rho}{\to} F\to F/F_+^\times\simeq\sign\,, \quad  \rho(T)= \alpha+t\epsilon
\end{equation}
for some $t\in \R$.
\begin{lem}\label{compree} Let $\varphi_1,\varphi_2\in\Hom(\Q[T],\sign)$ be two finite functions with $Re(\varphi_j)\in\R$, $j=1,2$. Then\vspace{.05in}

$(1)$~ $\varphi_1*_{\Delta^+}\varphi_2$
and $\varphi_1*_{\Delta^\times}\varphi_2$ are non-empty sets.\vspace{.05in}

$(2)$~$Re(\varphi) = Re(\varphi_1)+ Re(\varphi_2),\quad\forall\varphi\in\varphi_1*_{\Delta^+}\varphi_2$.\vspace{.05in}

$(3)$~ $Re(\varphi) = Re(\varphi_1)Re(\varphi_2),\quad \forall\varphi\in\varphi_1*_{\Delta^\times}\varphi_2$.
\end{lem}
\proof Lemmas \ref{tml} and \eqref{19k} ensure that $\varphi_1*_{\Delta^+}\varphi_2$
and $\varphi_1*_{\Delta^\times}\varphi_2$ are non-empty sets. For $j=1,2$, let $\alpha_j = Re(\varphi_j)$.  Let $\varphi\in\varphi_1*_{\Delta^+}\varphi_2$.  For $a_j\in \Q$, one has
$
\Delta^+(T-a_1-a_2) = (T-a_1)\otimes 1 + 1\otimes(T-a_2)
$
so that
\begin{equation}\label{13k}
\varphi(T-a_1-a_2)\in\varphi_1(T-a_1)\varphi_2(1)+\varphi_1(1)\varphi_2(T-a_2)
\end{equation}
By applying \eqref{cut}, the above inclusion shows that if $a_j<\alpha_j$ then $\varphi(T-a_1-a_2) = 1$ and $a_1+a_2 \le Re(\varphi)$. Similarly, if $a_j > \alpha_j$, one gets  that $\varphi(T-a_1-a_2) = -1$ and $a_1+a_2 \geq Re(\varphi)$. This shows that $Re(\varphi) = \alpha_1+\alpha_2$ and $(2)$ follows.

Let $\varphi\in\varphi_1*_{\Delta^\times}\varphi_2$.
For $a,b\in \Q$, one has
$$
\Delta^\times(T-ab) = T\otimes T - ab = (T-a)\otimes T + a\otimes(T-b).
$$
Thus
\begin{equation}\label{14k}
\varphi(T-ab) \in \varphi_1(T-a)\varphi_2(T) + \varphi_1(a)\varphi_2(T-b).
\end{equation}
To prove $(3)$, we first assume that $\alpha_j = Re(\varphi_j) >0$. This implies  $\varphi_j(T)=1$.
Thus, for $0<a<\alpha_1$, $0<b<\alpha_2$, one gets
$\varphi(T-ab) \in \varphi_1(T-a)\varphi_2(T) + \varphi_1(a)\varphi_2(T-b) = 1\cdot 1 + 1\cdot 1 = 1$.
Similarly, for $\alpha_1<a$, $\alpha_2<b$, one gets in the same way $\varphi(T-ab) = -1$. Thus, one derives $Re(\varphi) = \alpha_1\alpha_2$ when $\alpha_j>0$. One can check easily that the same result holds when
 $\alpha_j\neq 0$. Indeed, one can change the signs using  the automorphism $\sigma\in\Aut(\Q[T])$, $\sigma(T)= -T$. This way one obtains
 \begin{equation}\label{17k}
Re(\varphi^\sigma) = -Re(\varphi)
\end{equation}
where $\varphi^\sigma = \varphi\circ\sigma$, for $\varphi\in\Hom(\Q[T],\sign)$. Moreover it follows from
\begin{equation}\label{16k}
\Delta^\times\circ\sigma = (\sigma\otimes id)\circ\Delta^\times = (id \otimes \sigma)\circ\Delta^\times
\end{equation}
that
\begin{equation}\label{15k}
\varphi_1^\sigma*_{\Delta^\times}\varphi_2 = (\varphi_1*_{\Delta^\times}\varphi_2)^\sigma.
\end{equation}

We still have to consider the case when one of the $\alpha_j$'s vanishes. Assume $\alpha_1=0$. If $\varphi_1(T) = 0$, then for any $\varphi\in\varphi_1*_{\Delta^\times}\varphi_2$ one has $\varphi(T)=0$ since $\Delta^\times(T) = T\otimes T$. This shows that $0=\omega_0$ is an absorbing element. Then, by using the above change of signs automorphism (if needed), it is enough to consider the case $\varphi_1=\omega_0^+$ and $\varphi_2(T)=1$. One then obtains $\alpha_2\geq 0$. Let $b\in \Q$, $b> \alpha_2$ and $c>0$. Let $a>0$, $ab <c$.   Then, by applying \eqref{14k} we conclude
\[
\varphi(T-ab)\in \varphi_1(T-a)\varphi_2(T)+\varphi_1(a)\varphi_2(T-b)=-1.
\]
Thus, one gets $\varphi(T-c)=-1$. By  \eqref{cut} this gives $c\geq Re(\varphi)$ and since $c>0$ is arbitrary, one gets $Re(\varphi)\le 0$. But since $\varphi_2(T)=1$  one obtains $\varphi(T)\in\varphi_1(T)\varphi_2(T)=1$  for any $\varphi\in\varphi_1*_{\Delta^\times}\varphi_2$. This implies  $Re(\varphi)\ge 0$ and one then derives $Re(\varphi)=0$. Then, $(3)$ follows.
\endproof

\begin{lem}\label{transop} Let  $\varphi_1,\varphi_2\in\Hom(\Q[T],\sign)$, with $\alpha_j=Re(\varphi_j)\in\R$ ($j=1,2$) and assume $\alpha_1\notin\bar \Q$. Then\vspace{.05in}

$(1)$~$\varphi_1*_{\Delta^+}\varphi_2=Re^{-1}(\alpha_1+\alpha_2)$.\vspace{.05in}

$(2)$~If $\alpha_2\neq 0$, one has $\varphi_1*_{\Delta^\times}\varphi_2=Re^{-1}(\alpha_1\alpha_2)$.\vspace{.05in}

$(3)$~If $\alpha_2=0$, one has $\varphi_1*_{\Delta^\times}\varphi_2=\omega_0^s$ where $s=\ss(\alpha_1) \varphi_2(T)$.
\end{lem}

\proof Lemma \ref{compree} determines $\varphi_1*_{\Delta^+}\varphi_2$ (\resp $\varphi_1*_{\Delta^\times}\varphi_2$) if $\alpha_1+\alpha_2\notin\bar\Q$ (\resp $\alpha_1\alpha_2\notin \bar\Q$), since  in that case  $Re^{-1}(\alpha_1+\alpha_2)$
(resp. $Re^{-1}(\alpha_1\alpha_2)$) is reduced to a single element. To show (1) we can thus  assume that $\alpha_1+\alpha_2\in\bar\Q$. Then $\alpha_j\notin \bar\Q$ and $\varphi_j$ can be written in the form \eqref{19k} for any choice of $t$. One can then use Lemma~\ref{tml} to show that the three $\varphi$'s with $Re(\varphi)=\alpha_1+\alpha_2$ are in $\varphi_1*_{\Delta^+}\varphi_2$. A similar argument applies to show (2), since one can choose $t_j\in \R$ such that $p(T)\mapsto \ss(p((\alpha_1+t_1\epsilon)(\alpha_2+t_2\epsilon))$ yields any given element of $Re^{-1}(\alpha_1\alpha_2)$.

Finally, to prove (3) we note that if $\alpha_2=0$ then $\varphi\in\varphi_1*_{\Delta^\times}\varphi_2$ is determined by $\varphi(T)\in \sign$ and this
is given by $\varphi_1(T)\varphi_2(T)=\ss(\alpha_1) \varphi_2(T)$. \endproof

Lemma \ref{transop} determines the hyperoperations except when both
 $\alpha_j\in\bar\Q$. To discuss thoroughly the case $\alpha_j\in\bar\Q$, we first need to state the following two lemmas.

 \begin{lem}\label{posdec} Let $A(X,Y)\in\Q[X,Y]$ and let $\alpha_j\in\R$ such that $A(\alpha_1,\alpha_2)>0$. Then, there exist  polynomials $L_j\in\Q[X]$ and $R_j\in\Q[Y]$ such that
\begin{equation}\label{decplus}
A(X,Y) = \sum L_j(X)R_j(Y),~L_j(\alpha_1)>0,~R_j(\alpha_2)>0.
\end{equation}
\end{lem}
\proof One has
\[
A(X,Y)=\sum a_{n,m}X^nY^m,\quad\sum a_{n,m}\alpha_1^n\alpha_2^m>0.
\]
We take $\epsilon_{n,m}>0$ such that $\xi = \sum a_{n,m}\alpha_1^n\alpha_2^m-\sum\epsilon_{n,m}>0$ and write
\[
A(X,Y) = \sum(a_{n,m}X^nY^m-a_{n,m}\alpha_1^n\alpha_2^m+\epsilon_{n,m})+\xi.
\]
This shows that it is enough to prove the lemma in the special case  $$A(X,Y)=aX^nY^m-b \,, \ a\alpha_1^n\alpha_2^m>b.$$ In this case one can find $\lambda,\mu \in \Q$ such that
$$
a(\alpha_1^n-\lambda)>0\,, \ (\alpha_2^m-\mu)>0\,, \ A(\alpha_1,\alpha_2)>a(\alpha_1^n-\lambda)(\alpha_2^m-\mu).
$$
Indeed, this follows by using $\lambda=\alpha_1^n-\epsilon_1$, $\mu=\alpha_2^m-\epsilon_2$ and taking $\epsilon_j$ with $a\epsilon_1>0$ and $\epsilon_2>0$ small enough so that
 $$
 a\epsilon_1\epsilon_2<A(\alpha_1,\alpha_2).
 $$
  One has the decomposition \eqref{decplus}
  $$
  a(X^n-\lambda)(Y^m-\mu)=L_1(X)R_1(Y) \,, \ L_1(\alpha_1)>0\,, \ R_1(\alpha_2)>0\,.
  $$
  It is thus enough to look at the remainder
 \[
A_1(X,Y)=A(X,Y) - a(X^n-\lambda)(Y^m-\mu).
\]
Then, $A_1(X,Y)$ is a polynomial of the form $$A_1(X,Y)=\alpha X^n+\beta Y^m +\delta\,, \ A_1(\alpha_1,\alpha_2)=\alpha\alpha_1^n+\beta\alpha_2^m+\delta>0.$$ Thus, one derives the description
$$
A_1(X,Y)=(\alpha(X^n-\alpha_1^n)+\delta_1)+(\beta(Y^m-\alpha_2^m)+\delta_2)\,, \ \delta_j>0.
$$
which is of the form \eqref{decplus} since $$A_1(X,Y)=L_2(X) +   R_2(Y)\,,
L_2(\alpha_1)>0\,, \ R_2(\alpha_2)>0.$$
\endproof

\begin{lem} \label{ideal} Let $p_1$ and $p_2$ be two irreducible polynomials in $\Q[T]$. Then a polynomial $p(X,Y)$ belongs to the ideal
\begin{equation}\label{32k}
J=\{p_1(X)A(X,Y)+p_2(Y)B(X,Y)\mid A,B\in \Q[X,Y]\}\subset\Q[X,Y]
\end{equation}
if and only if $p(\alpha,\beta)=0$, $\forall~\alpha\in Z(p_1)$, $\beta\in Z(p_2)$.
\end{lem}
\proof By the Nullstellensatz, it is enough to show that $\Q[X,Y]/J$ is a reduced algebra. By construction  $\Q[X,Y]/J\simeq\Q[X]/(p_1)\otimes\Q[Y]/(p_2)$ and since $p_j$ is irreducible, this is the tensor product of two fields, hence it is reduced.
\endproof

\begin{prop}\label{hypsign}The hyperaddition $\star_{\Delta^{+}}$ on  functions $\varphi\in \Hom(\Q[T],\sign)$ satisfying the condition $Re(\varphi)\in \bar\Q\cap \R$ coincide with the hyperaddition in $\sign \times (\bar\Q\cap \R)$.
\end{prop}

\proof
 We first consider the case $\varphi_j=\omega_{\alpha_j}$. Then,  by applying Lemma~\ref{sml} one gets that $\omega_{\alpha_1+\alpha_2}$ is the only element  in $\varphi_1*_{\Delta^{+}}\varphi_2$, since its kernel  is non trivial. Using \eqref{19k} one gets, with $\alpha = \alpha_1+\alpha_2$
\begin{equation}\label{threecases}
\omega^+_{\alpha_1}+\omega^-_{\alpha_2} = \{\omega_\alpha,\omega^+_\alpha,\omega^-_\alpha\},
\ \
\omega^+_\alpha\in\omega_{\alpha_1}+\omega^+_{\alpha_2},\quad \omega^+_{\alpha}\in\omega^+_{\alpha_1}+\omega^+_{\alpha_2}.
\end{equation}
We need to show that there are no other solutions in the last two cases. Let $p_j(T)\in\Q[T]$ ($j=1,2$) be two irreducible polynomials with $p_j(\alpha_j)=0$ and $p'_j(\alpha_j)>0$. In general, the polynomial $q(Z)$ obtained in \eqref{resu}, \ie as the resultant in the variable $T$ of $p_1(T)$ and $p_2(Z-T)$, may have $\alpha_1+\alpha_2$ as a multiple root. We replace $q$ by the product $p=\prod q_j$ of the irreducible factors $q_j$ which appear in the decomposition  $q=\prod q_j^{n_j}$
of $q$ as a product of powers of prime factors. By construction $\alpha_1+\alpha_2$ is a simple root of $p$ and by multiplying by a non-zero scalar (if necessary) we can assume
\begin{equation}\label{22k}
p'(\alpha_1+\alpha_2)>0.
\end{equation}
By applying Lemma \ref{ideal} we deduce a decomposition
\begin{equation}\label{23k}
p(X+Y) = p_1(X)A(X,Y)+p_2(Y)B(X,Y).
\end{equation}
We want to show in this generality that
\begin{equation}\label{24k}
A(\alpha_1,\alpha_2)>0,\quad B(\alpha_1,\alpha_2)>0.
\end{equation}
We can differentiate \eqref{23k} with respect to $X$ and get
\[
p'(X+Y) = p_1'(X)A(X,Y)+p_1(X)\partial_XA(X,Y) + p_2(Y) \partial_XB(X,Y)
\]
and then take $X=\alpha_1$ and $Y=\alpha_2$ to get
\[
p'(\alpha_1+\alpha_2) = p_1'(\alpha_1)A(\alpha_1,\alpha_2).
\]
Thus, since $p'(\alpha_1+\alpha_2)>0$ and $p'_j(\alpha_j)>0$, we obtain $A(\alpha_1,\alpha_2)>0$. The same argument applies to $B(\alpha_1,\alpha_2)$ using $\partial_Y$.
Lemmas \ref{posdec}, \eqref{24k} and \eqref{23k} thus give a decomposition
\begin{equation}\label{posdecbis}
p(X+Y) = p_1(X)\sum A_j(X)B_j(Y)+p_2(Y)\sum C_j(X)D_j(Y)
\end{equation}
with
$$
 A_j(\alpha_1)>0,\ B_j(\alpha_2)>0,\ C_j(\alpha_1)>0,\ D_j(\alpha_2)>0.
 $$
This implies, using $Re(\varphi_j)=\alpha_j$, that
$$
 \varphi_1(A_j)=1,\ \varphi_2(B_j)=1,\ \varphi_1(C_j)=1,\ \varphi_2(D_j)=1.
 $$
By using \eqref{posdecbis} we get
$$
\varphi(p)\in \varphi_1(p_1)+\varphi_2(p_2)\qqq \varphi \in \varphi_1*_{\Delta^+}\varphi_2
$$
and by Corollary \ref{sign} (4), one concludes
$$
\omega_{\alpha_1}+\omega^+_{\alpha_2}=\omega^+_\alpha,\quad \omega^+_{\alpha_1}+\omega^+_{\alpha_2}=\omega^+_{\alpha}.
$$
This completes the table of hyperaddition.
 \endproof

\begin{prop}\label{hypsignmult}The hypermultiplication $\star_{\Delta^\times}$ on  functions $\varphi\in \Hom(\Q[T],\sign)$ such that $Re(\varphi)\in \bar\Q\cap \R$ is given, for $\alpha_j\neq 0$ ($j=1,2$), by
\begin{equation}\label{hypprod}
  \omega^{s_1}_{\alpha_1}\star_{\Delta^\times}\omega^{s_2}_{\alpha_2}=\{\omega^{s}_{\alpha_1\alpha_2}
  \mid s\in \ss(\alpha_2) s_1+\ss(\alpha_1)s_2\subset \sign\}\qqq s_j\in \sign.
\end{equation}
For $\alpha_1=0$ and any value of $\alpha_2\in \bar\Q\cap \R$, one has
\begin{equation}\label{hypprod0}
\omega^{s_1}_{0}\star_{\Delta^\times}\omega^{s_2}_{\alpha_2}=\omega^{s}_{0}\,, \ s=s_1\omega^{s_2}_{\alpha_2}(T).
\end{equation}
\end{prop}
\proof  By representing $\omega^{s_j}_{\alpha_j}$ in the form \eqref{19k}, one gets that all elements of the right hand side of \eqref{hypprod} belong to the left hand side. To obtain the other inclusion, we first assume that $\alpha_j>0$. We proceed as in the proof of Proposition \ref{hypsign} and obtain a polynomial with simple roots which admits as roots all the products $g_1(\alpha_1)g_2(\alpha_2)$ of the conjugates of the $\alpha_j$. By using the same notation as before, one gets a decomposition
\begin{equation}\label{23kbis}
p(XY) = p_1(X)A(X,Y)+p_2(Y)B(X,Y).
\end{equation}
Since $\alpha_j>0$, the same proof using differentiation shows that \eqref{24k} holds.
One then gets
$$
\varphi(p)\in \varphi_1(p_1)+\varphi_2(p_2)\qqq \varphi \in \varphi_1*_{\Delta^\times}\varphi_2
$$
which shows the required conclusion. By using \eqref{15k}, one then obtains the general case, when $\alpha_j\neq 0$ have arbitrary signs.
Finally \eqref{hypprod0} follows from  $\Delta^\times(T)=T\otimes T$.
\endproof

Incidentally, we note that the product $A_1\times A_2$ of two canonical hypergroups $(A_1,+_1)$ and $(A_2,+_2)$ endowed with the hyperaddition
\begin{equation}\label{prodhyp}
(\xi_1,\xi_2)+(\eta_1,\eta_2)  = \{(\alpha,\beta), \alpha\in \xi_1+_1\eta_1,~\beta\in \xi_2+_2\eta_2\}
\end{equation}
is a canonical hypergroup.

\begin{lem}\label{glue} Let $A$ be an abelian group, and $B\subset A$ be a subgroup. Let $C$ be the quotient of the product $A\times \sign$ by the
equivalence relation
\begin{equation}\label{29k}
(\alpha,s)~\sim~(\alpha,0)\qquad\forall\alpha\notin B.
\end{equation}
We let $\epsilon: A\times \sign \to C$ be the quotient map and endow $A\times \sign$
with the hyperaddition \eqref{prodhyp}. Then, with the following hyperlaw,
$C$ is a canonical hypergroup:
\begin{equation}\label{30k}
x+y = \{\epsilon(\xi+\eta),~\epsilon(\xi)=x, \epsilon(\eta)=y\}.
\end{equation}
\end{lem}
\proof
For $x=(\alpha,s)$ and $y = (\beta,t)$, the only case when one needs to take the union on representatives $\xi$ with $\epsilon(\xi)=x$ and $\eta$ with $\epsilon(\eta)=y$ is when both $\alpha,\beta\notin B$ but $\alpha+\beta\in B$. Note also that for the equivalence relation \eqref{29k} the subset
\begin{equation}\label{31k}
\cup(\xi+\eta),\quad\epsilon(\xi)=x,~\epsilon(\eta)=y
\end{equation}
is saturated. Indeed, if $\alpha+\beta\in B$ there is nothing to prove since the equivalence is trivial. Thus one can assume $\alpha+\beta\notin B$ and  say $\alpha\notin B$. Then, one can choose $\xi = (\alpha,s)$ with all values $s\in\sign$ so that the whole fiber above $\alpha+\beta$ appears in $\xi+\eta$.

This produces the associativity of the hyperaddition in  $C=(A\times\sign)/\sim$ (\cf\eqref{29k}) by using the associativity in $A\times \sign$. The uniqueness of the additive inverse follows from $\epsilon(\xi)=0\implies\xi =0$. Finally, the reversibility property follows from \eqref{31k} by using \eqref{29k}. \endproof

We denote by $A\times_{B^c}\sign$ the canonical hypergroup obtained from Lemma \ref{glue}.

\begin{thm} The functions $\cD_{\rm finite}(\sign)$ form under hyperaddition $\star_{\Delta^{+}}$ a canonical hypergroup isomorphic to $\R\times_{\bar\Q^c}\sign$.

The subset $\cD_{\rm finite}^{\times}(\sign)=\{\varphi\in\cD_{\rm finite}(\sign)|Re(\varphi)\neq 0\}$, with the hypermultiplication $\star_{\Delta^{\times}}$ is a canonical hypergroup isomorphic to $\R^\times\times_{\bar\Q^c}\sign$.
\end{thm}
\proof
The first statement is obtained using Lemma \ref{transop} (1) and Proposition \ref{hypsign}.
To show the second statement, we introduce the map
\begin{equation}\label{sigmamap}
 \sigma:\cD_{\rm finite}^{\times}(\sign)\to\R^\times\times_{\bar\Q^c}\sign,\quad   \sigma(\omega_\alpha^s)=(\alpha,\ss(\alpha)s).
\end{equation}
Then,  Lemma \ref{transop} (2) and \eqref{hypprod} show that $\sigma$ is an isomorphism.
\endproof

\begin{rem}\label{ssub}{\rm
The set $\cD_{\rm finite}^{0}(\sign)=\{\varphi\in\cD_{\rm finite}(\sign)|Re(\varphi)= 0\}$ forms
an ideal of the hyperstructure $\cD_{\rm finite}(\sign)$ with $\cD_{\rm finite}^{0}(\sign)\simeq \sign$.
 The quotient hyperstructure is isomorphic to the field of real numbers:
 \begin{equation}\label{rrr}
   \cD_{\rm finite} (\sign)/\sign\simeq\R.
 \end{equation}
}\end{rem}

\section{The hyperring  of ad\`ele classes and its arithmetic} \label{hyper}

The theory of hyperrings allows one to understand the additive structure of the multiplicative monoid $\A_\K/\K^\times$ of  ad\`ele classes of a global field $\K$ and to obtain a new algebraic understanding of the ad\`ele class space as a $\kras$-algebra. Indeed, by Theorem \ref{wagnad} and Corollary \ref{wagnad1},  the quotient $\ads=\A_\K/\K^\times$ of the commutative ring $\A_\K$ by the action by multiplication of $\K^\times$, is a hyperring  and the Krasner hyperfield $\kras$ is embedded  in $\ads$. In short,   $\ads$ is the $\kras$-algebra $\ads=\A_\K\otimes_\K\kras$, obtained by extension of scalars using the unique homomorphism $\K\to \kras$.

In this section we shall review the most important arithmetic properties of the hyperring $\H_\K$.  The set $P(\ads)$ of the {\em prime} elements of the hyperring   $\ads$ inherits a natural structure of groupoid with the product given by multiplication and  units the set of places of $\K$. The id\`ele class group $C_\K=\ads^\times$  acts by multiplication on $P(\ads)$.
For a global field of positive characteristic, the action of the  units
$\H_\K^\times$ on the {\em prime elements} of $\H_\K$ corresponds, by class-field theory, to the action of
the abelianized  Weil group  on
the space of valuations of the maximal abelian
extension of $\K$, \ie on the space of the (closed) points of the
corresponding projective tower of algebraic curves. This construction determines the maximal abelian cover of the projective algebraic curve  with function field $\K$. Then, the sub-groupoid of loops of the fundamental groupoid associated to the afore mentioned tower   is equivariantly isomorphic to $P(\ads)$.

When ${\rm char}(\K)=0$, the above geometric interpretation is no longer available. On the other hand, the arithmetic of the hyperring $\ads$ continues to hold and the  groupoid $P(\ads)$  appears to be a natural substitute for the above groupoid of loops  since it  also supports an interpretation of the explicit formulas of Riemann-Weil.

\subsection{The space of closed prime ideals of $\ads$}

Let $\K$ be a global field and $\Sigma(\K)$ the set of places of $\K$. The one-to-one correspondence between  subsets $Z\subset \Sigma(\K)$ and closed ideals of $\A_\K$ (for the locally compact topology)  given by
\begin{equation}\label{idealsads}
Z\mapsto J_Z=\{x=(x_v)\in \A_\K\mid x_w=0\qqq w\in Z\}.
\end{equation}

determines, when $Z = \{w\}$ ($w\in\Sigma(\K)$), a one-to-one relation  between the places of $\K$ and the {\em prime, closed} ideals of the hyperring $\H_\K=\A_\K/\K^\times$
\begin{equation}\label{prideal}
\Sigma(\K)\ni w~\leftrightarrow~\ffp_w=\{x\in \H_\K\,|\, x_w=0\}\subset\Sp(\H_\K).
\end{equation}
Notice that the ideal $\ffp_w$  is well defined since  the condition for an ad\`ele to vanish at a place is invariant under multiplication by elements in $\K^\times$.

The additive structure of $\ads$ plays a key role in the above relation since when viewed as a multiplicative monoid, the ad\`ele class space $\A_\K/\K^\times$ has many more prime ideals than when it is viewed  as a $\kras$-algebra. In fact, in a monoid  {\em any} union of prime ideals is still a prime ideal and this fact implies that all subsets of the set of places determine a prime ideal of the monoid of ad\`ele classes.

\subsection{The groupoid of  prime elements of $\ads$}\label{space}

 In a  hyperring $R$,  an element $a\in R$ is said to be {\em prime} if the ideal $aR$ is a prime ideal.
Let $P(\ads)$ be the set of prime elements of the hyperring $\ads=\A_\K/\K^\times$. Each prime element $a\in P(\ads)$ determines a principal, prime ideal $\ffp=a\ads\subset \H_\K$.  The following result establishes a precise description of  such ideals of $\H_\K$. We refer to \cite{wagner} Theorem~7.9 for further details.

\begin{thm}\label{classprime}~$1)$~Any  principal prime ideal $J = a\H_\K$ of $\ads$ is equal to $\ffp_w$ for a  unique place $w\in\Sigma(\K)$.\vspace{.05in}

$2)$~The group $C_\K=\A_\K^\times/\K^\times$ acts transitively on the generators of the principal prime ideal $\ffp_w$.\vspace{.05in}

$3)$~The isotropy subgroup of any generator of the  prime ideal $\ffp_w$ is $\K_w^\times\subset C_\K$.
\end{thm}

Let  $s:P(\ads)\to \Sigma_\K$  be the map that associates to a prime element $a\in \ads$ the unique place $w$ such that $\ffp_w=a\ads$. Then $P(\ads)$ with range and source maps equal to $s$ and partial product given by the product in the hyperring $\ads$,
is a {\em groupoid}. The product of two prime elements is  a prime element only when the two factors generate the same ideal, \ie sit over the same place. Moreover over each place $v$ there exists a unique  idempotent $p_v\in P(\ads)$ (\ie $p^2_v=p_v$).

\subsection{The groupoid of loops and $P(\ads)$ in characteristic $p\neq 0$}

Let $\K$ be  a global field of characteristic $p> 0$ \ie a function field over a constant field $\F_q\subset \K$. We fix a separable closure $\bar \K$ of $\K$ and let
 $\K^{\rm ab}\subset \bar \K$ be the maximal abelian extension of $\K$. Let $\bar \F_q\subset \bar \K$ be the algebraic closure of  $\F_q$. We denote by $\cW^{\rm ab}\subset {\rm Gal}( \K^{\rm ab}:\K)$  the abelianized  Weil group, \ie the subgroup of elements of ${\rm Gal}( \K^{\rm ab}:\K)$ whose restriction to $\bar \F_q$ is an integral power of the Frobenius.

Let ${\rm Val}(\K^{\rm ab})$ be the space of all  valuations of $\K^{\rm ab}$. By restriction to $\K\subset \K^{\rm ab}$ one obtains a natural map
\begin{equation}\label{mapp}
    p\;:\; {\rm Val}(\K^{\rm ab})\to \Sigma(\K)\,, \ \ p(v)=v|_\K.
\end{equation}
By construction, the action of ${\rm Gal}( \K^{\rm ab}:\K)$ on ${\rm Val}(\K^{\rm ab})$ preserves the map $p$.

Let $w\in \Sigma(\K)$. Then, it follows from standard results of class-field theory that the abelianized  Weil group $\cW^{\rm ab}$ acts {\em transitively} on the fiber $p^{-1}(w)$ of $p$ and that the isotropy subgroup of an element in the fiber $p^{-1}(w)$ coincides with the abelianized  local Weil   group $\cW^{\rm ab}_w\subset \cW^{\rm ab}$.\vspace{.02in}

We now implement the geometric language.  Given an extension $\bar \F_q\subset E$  of transcendence degree $1$, it is a well-known fact that the  space ${\rm Val}(E)$  of valuations of $E$ coincides with the set of (closed) points of the unique projective nonsingular algebraic curve with function field $E$. Moreover, one also knows  (\cf \cite{Hart} Corollary 6.12) that the category of nonsingular projective algebraic curves and dominant morphisms is equivalent to the category of function fields of dimension one over $\bar \F_q$.  One also knows that the maximal abelian extension $\K^{\rm ab}$ of $\K$ is an inductive limit of extensions $E$  of $\bar \F_q$ of transcendence degree $1$.
Thus the space ${\rm Val}(\K^{\rm ab})$ of valuations of $\K^{\rm ab}$, endowed with the action of the abelianized  Weil group $\cW^{\rm ab}\subset{\rm Gal}(\K^{\rm ab}:\K)$, inherits the structure of a projective limit of projective nonsingular curves. This construction determines the maximal abelian cover $\pi:X^{\rm ab}\to X$ of the non singular projective curve $X$ over $\F_q$ with function field  $\K$.\vspace{.05in}

Let $\pi:\tilde X\to X$ be a Galois covering of $X$ with Galois group $W$. The {\em fundamental groupoid} of $\pi$ is by definition the quotient $\Pi_1=(\tilde X\times \tilde X)/W$ of $\tilde X\times \tilde X$ by the diagonal action of $W$ on the self-product. The (canonical) range and source maps: $r$ and $s$ are defined by the two projections
\begin{equation}\label{rs}
    r(\tilde x,\tilde y)=x\,, \ s(\tilde x,\tilde y)=y.
\end{equation}
Let us consider the subgroupoid of {\em loops} \ie
\begin{equation}\label{sub}
    \Pi_1'=\{\gamma\in \Pi_1\mid r(\gamma)=s(\gamma)\}.
\end{equation}
 Each fiber of the natural projection $r=s:\Pi_1'\to X$ is a group. Moreover, if $W$ is an abelian group  one defines the following natural action of $W$ on $\Pi_1'$
\begin{equation}\label{actionW}
    w\cdot (\tilde x,\tilde y)=(w\tilde x,\tilde y)=(\tilde x,w^{-1}\tilde y).
\end{equation}
We consider the maximal abelian cover $\pi:X^{\rm ab}\to X$ of the non singular projective curve $X$ over $\F_q$ with function field $\K$.
We view $X$ as a scheme over $\F_q$. In this case, we let $W=\cW^{\rm ab}\subset{\rm Gal}(\K^{\rm ab}:\K)$  be the abelianized Weil group. Even though the  maximal abelian cover $\pi:X^{\rm ab}\to X$ is ramified, its  loop groupoid $\Pi_1^{\rm ab}(X)'$ continues to make sense. Since the two projections from $X^{\rm ab}\times X^{\rm ab}$ to $X$ are $W$-invariant, $\Pi_1^{\rm ab}(X)'$ is the quotient of the fibered product $X^{\rm ab}\times_X X^{\rm ab}$ by the diagonal action of $W$. We identify the closed points of $X^{\rm ab}\times_X X^{\rm ab}$ with pairs of valuations of $\K^{\rm ab}$ whose restrictions to $\K$ are the same.
The following refinement of Proposition 8.13 of
\cite{CCM2} holds

\begin{thm}\label{ccm2prop} Let $\K$ be a  global field of characteristic $p\neq 0$, and let $X$ be the corresponding non-singular projective algebraic curve over $\F_q$.\vspace{.05in}

 $\bullet$~The loop groupoid $\Pi_1^{\rm ab}(X)'$ is canonically isomorphic to the groupoid $P(\ads)$ of prime elements of the hyperring $\ads=\A_\K/\K^\times$.\vspace{.05in}

    $\bullet$~The above isomorphism $\Pi_1^{\rm ab}(X)'\simeq P(\ads)$  is equivariant for the action of $W$ on $\Pi_1^{\rm ab}(X)'$  and the action of the units $\ads^\times=C_\K$ on prime elements by multiplication.
 \end{thm}


\begin{thebibliography}{99}


\bibitem{Beutel} A.~Beutelspacher,
{\em Projective planes}.
Handbook of incidence geometry, 107--136, North-Holland, Amsterdam, 1995.

\bibitem{Bourbaki} N.~Bourbaki, {\em  Algebra II. Chapters 4--7}. Translated from the 1981 French edition by P. M. Cohn and J. Howie. Reprint of the 1990 English edition. Elements of Mathematics (Berlin). Springer-Verlag, Berlin, 2003.

    \bibitem{Cartier} P.~Cartier, {\em Analyse num\'erique d'un probl\`eme de valeurs propres a haute pr\'ecision, applications aux fonctions automorphes}, Preprint IHES, (1978).


\bibitem{Co-zeta} A.~Connes, {\em Trace formula in noncommutative geometry and the zeros of the Riemann zeta function},  Selecta Math. (N.S.)  5  (1999),  no. 1, 29--106.


\bibitem{ak} A.~Connes, C.~Consani,  {\em On the notion of geometry over $\F_1$}, to appear in Journal of Algebraic Geometry; arXiv08092926v2 [mathAG].

\bibitem{announc3} A.~Connes, C.~Consani, {\em Schemes over $\F_1$ and zeta functions}, to appear in Compositio Mathematica; arXiv:0903.2024v3 [mathAG,NT].

\bibitem{jamifine} A.~Connes, C.~Consani, {\em Characteristic $1$, entropy and the absolute point}, to appear in the Proceedings of the 21st JAMI Conference, Baltimore 2009, JHUP; arXiv:0911.3537v1 [mathAG].

\bibitem{wagner} A.~Connes, C.~Consani, {\em The hyperring of ad\`ele classes}, to appear in Journal of Number Theory; arXiv:1001.4260 [mathAG,NT].


\bibitem{CCM} A.~Connes, C.~Consani, M.~Marcolli, {\em
Noncommutative geometry and motives: the thermodynamics of endomotives},
Advances in Math. 214 (2) (2007), 761--831.

\bibitem{CCM2} A.~Connes, C.~Consani, M.~Marcolli, {\em The Weil proof
and the geometry of the adeles class space}, in ``Algebra, Arithmetic
and Geometry -- Manin Festschrift'', Progress in Mathematics, Birkh\"auser
(2010).

\bibitem{ccm}  A.~Connes, C.~Consani, M.~Marcolli, {\em
Fun with $\F_1$}, Journal of Number Theory 129 (2009) 1532--1561.

\bibitem{CMbook} A.~Connes, M.~Marcolli,  {\em  Noncommutative Geometry, Quantum Fields, and Motives},
Colloquium Publications, Vol.55, American Mathematical Society, 2008.






 \bibitem{deit} A.~Deitmar, {\em Schemes over F1}, in Number Fields and Function Fields--Two Parallel Worlds. Ed. by G. van der Geer, B. Moonen, R. Schoof. Progr. in
Math, vol. 239, 2005.

\bibitem{demgab} M.~Demazure, P.~Gabriel, {\em  Groupes alg\'ebriques}, Masson \& CIE, \'Editeur Paris 1970.

\bibitem{Ellers} E.~Ellers, H.~Karzel,
{\em Involutorische Geometrien}. (German)
Abh. Math. Sem. Univ. Hamburg 25 1961 93--104.



\bibitem{goss} D. Goss, {\em Basic Structures of Function Field Arithmetic}. Ergebnisse der Mathematik und ihrer Grenzgebiete (3) [Results in Mathematics and Related Areas (3)], 35. Springer-Verlag, Berlin, 1996. xiv+422 pp.

    \bibitem{gui}  V. Guillemin,  {\em Lectures on spectral theory of
elliptic operators},   Duke Math. J.,  Vol. 44, 3 (1977), 485--517.



\bibitem{Hall}  M.~Hall, {\em
Cyclic projective planes}. Duke Math. J. 14, (1947), 1079--1090

\bibitem{Hart} R.~Hartshorne  {\em Algebraic Geometry}, Graduate Texts in Mathematics 52, Springer-Verlag, New York Heidelberg Berlin 1977.


\bibitem{Ingham} A.~Ingham, {\em The distribution of prime numbers} With a foreword by R. C. Vaughan. Cambridge Mathematical Library. Cambridge University Press, Cambridge, 1990. xx+114 pp.




\bibitem{Karzel} H.~Karzel, {\em Bericht \"{u}ber projektive Inzidenzgruppen}.
(German) Jber. Deutsch. Math.-Verein. 67 1964/1965 Abt. 1, 58--92.


\bibitem{Kato} K.~Kato, {\em Toric Singularities}
 American Journal of Mathematics, Vol. 116, No. 5 (Oct., 1994), 1073--1099.




\bibitem{Krasner} M.~Krasner, {\em Approximation des corps valu\'es complets de caract\'eristique
$p\not=0$ par ceux de caract\'eristique $0$},  (French) 1957 Colloque
d'alg\`ebre sup\'erieure, tenu \`a Bruxelles du 19 au 22 d\'ecembre 1956 pp.
129--206 Centre Belge de Recherches Math\'ematiques \'Etablissements
Ceuterick, Louvain; Librairie Gauthier-Villars, Paris.


\bibitem{Krasner1} M.~Krasner, {\em A class of hyperrings and hyperfields}. Internat. J. Math. Math. Sci. 6 (1983), no. 2, 307--311.

    \bibitem{Ku} N.~Kurokawa, {\em Multiple zeta functions: an example.}
 in Zeta functions in geometry (Tokyo, 1990), Adv. Stud. Pure Math., 21, Kinokuniya, Tokyo, 1992, pp.  219--226.


\bibitem{KOW} N.~Kurokawa,  H.~Ochiai, A.~Wakayama, {\em Absolute Derivations
and Zeta Functions} Documenta Math. Extra Volume: Kazuya Kato's
Fiftieth Birthday (2003) 565--584.





\bibitem{Manin} Y.~I.~Manin, {\em Lectures on zeta functions and motives (according to Deninger and Kurokawa)} Columbia University Number-Theory Seminar (1992),
Ast\'erisque No.~228 (1995), 4, 121--163.


\bibitem{Meyer} R.~Meyer, {\em On a representation of the idele class group
related to primes and zeros of $L$-functions}.  Duke Math. J. Vol.127 (2005),
N.3, 519--595.








\bibitem{Prenowitz} W.~Prenowitz,  {\em Projective Geometries as Multigroups}
 American Journal of Mathematics, Vol. 65, No. 2 (1943), pp. 235--256.






\bibitem{Soule} C.~Soul\'e, {\em Les vari\'et\'es sur le corps \`a un
\'el\'ement}. Mosc. Math. J. 4 (2004), no. 1, 217--244.






\bibitem{Thas} K.~Thas, D.~Zagier, {\em
Finite projective planes, Fermat curves, and Gaussian periods}.
J. Eur. Math. Soc. (JEMS) 10 (2008), no. 1, 173--190.

\bibitem{Tits} J.~Tits, {\em Sur les analogues alg\'ebriques des groupes semi-simples
complexes}. Colloque d'alg\`ebre sup\'erieure, tenu \`a Bruxelles du 19 au 22
d\'ecembre 1956, Centre Belge de Recherches Math\'ematiques \'Etablissements
Ceuterick, Louvain; Librairie Gauthier-Villars, Paris   (1957), 261--289.


\bibitem{TV} B.~T\"oen, M.~Vaqui\'e, {\em Au dessous de $\text{Spec}(\Z)$}. J. K-Theory 3 (3) (2009), 437--500.


\bibitem{Wagner} A.~Wagner, {\em
On perspectivities of finite projective planes}.
Math. Z 71 (1959), 113--123.

\bibitem{Wagner1} A.~Wagner, {\em
On collineation groups of projective spaces. I.}
Math. Z. 76 (1961), 411--426.


 \bibitem{Weil} A.~Weil, {\em Sur la th\'eorie du corps de classes}  J. math. Soc. Japan, t. 3, 1951,  1--35.



\end{thebibliography}
\end{document}